\documentclass[12pt,a4paper]{amsart}
\usepackage[usenames]{color}
\usepackage{times}
\usepackage{amsfonts}
\usepackage{amsthm}
\usepackage{amsmath}
\usepackage{psfrag}

\usepackage{times}
\usepackage{latexsym,color}
\usepackage{amssymb}
\usepackage{graphicx}

\usepackage{epsfig}

\usepackage{tikz}
\usepackage[all]{xy}
 \usepackage[hang]{caption}
 \usetikzlibrary{snakes}
\usetikzlibrary{arrows}

\advance \topmargin by -\headheight
\advance \topmargin by -\headsep
\evensidemargin \oddsidemargin
\usepackage{lineno}
\nolinenumbers

\newtheorem{theorem}{Theorem}
\newtheorem{lemma}[theorem]{Lemma}
\newtheorem{corollary}[theorem]{Corollary}
\newtheorem{proposition}[theorem]{Proposition}
\newtheorem{question}{Question}
\newtheorem{example}{Example}
\newtheorem{definition}{Definition}

\newtheorem{conj}{Conjecture}

\begin{document}

\title[AR clustering]{Clustering and Arnoux-Rauzy words
}

\author[S. Ferenczi]{S\'ebastien Ferenczi}
\address{Aix Marseille Universit\'e, CNRS, Centrale Marseille, Institut de Math\' ematiques de Marseille, I2M - UMR 7373\\13453 Marseille, France.}
\email{sebastien-simon.ferenczi@univ-amu.fr}

\author[L.Q. Zamboni]{Luca Q. Zamboni}
\address{Institut Camille Jordan,
Universit\'e Claude Bernard Lyon 1,
43 boulevard du 11 novembre 1918,
69622 Villeurbanne Cedex, France}
\email{zamboni@math.univ-lyon1.fr}

\subjclass[2010]{Primary 68R15}
\date{May 29, 2023}

\begin{abstract} We characterize  the clustering of a word under the Burrows-Wheeler transform
 in terms of the resolution of a bounded number of bispecial factors belonging to the language generated by all its powers. We use this criterion to compute,  in  every given Arnoux-Rauzy language on three letters, an explicit bound $K$ such that each word of length at least $K$ is not clustering; this bound is sharp for a set of Arnoux-Rauzy languages including the Tribonacci one. In the other direction, we characterize all standard Arnoux-Rauzy clustering words, and all perfectly clustering Arnoux-Rauzy words. We extend some results to episturmian languages, characterizing those which produce infinitely many clustering words, and to larger alphabets. 
 \end{abstract}

\maketitle

 In \cite{fz4}, the authors give a characterization of the {\it clustering} phenomenon for the {\it Burrows-Wheeler transform}, using a class of dynamical systems, the {\em interval exchange} transformations. This gives a way to build examples of clustering words, but is not very operative in deciding whether a given word clusters. Here, inspired by \cite{fhuz} but independently and  with purely combinatorial methods, we give, in Theorem \ref{car}, a characterization of the clustering of a primitive word $w$
 in terms of the resolution of a finite number of bispecial words of the language generated by all the $w^n$, $n>0$. 
 
 A very popular family of languages consists in the {\it Sturmian} languages; these are known since \cite{mrs} to be good producers of clustering words, though not all of their factors are clustering. As is  proved in \cite{fz4},  their natural generalizations, the  interval exchange languages, can also produce infinitely many clustering words.  In both cases, our Theorem \ref{car} gives a new criterion to identify those factors which cluster. Another well-known generalization of the Sturmian languages consists in the {\it Arnoux-Rauzy languages} on three letters, and the question of clustering of their factors was asked by Francesco Dolce at the Journées Montoises 2022:  in contrast with the previous cases, we are able to answer it by a broad negative, the Arnoux-Rauzy
 words are in general not clustering. 
 
 More precisely, for  every given Arnoux-Rauzy language, we compute, in Theorem \ref{arc}, an explicit bound $K$ such that each word of length at least $K$ is not clustering; by Corollary \ref{sha} this bound is indeed sharp for a family of Arnoux-Rauzy  languages including the Tribonacci language, but we can build  counter-examples for which this bound is not optimal. In the other direction, in each Arnoux-Rauzy language we want to find primitive clustering words $v$: it turns out that this is easier for {\it standard} Arnoux-Rauzy words, as we prove in  Proposition \ref{sq} that an Arnoux-Rauzy word $v$ is cyclically conjugate to a standard one if and only if $vv$ is an Arnoux-Rauzy word;  in  Proposition \ref{list}, and Corollary \ref{clist},  we are  able  to characterize those which are clustering. As a consequence, there  exist arbitrarily long primitive perfectly clustering Arnoux-Rauzy words, and every Arnoux-Rauzy language contains a primitive perfectly clustering word of length at least $22$. As for clustering Arnoux-Rauzy words not conjugate to standard ones, in Propositions \ref{ptb1} and \ref{ptb2} we characterize, using methods and results of \cite{sp}, those which cluster perfectly, but there are also infinitely many of them which cluster but not perfectly.
 
  Finally, we turn to generalizations of Arnoux-Rauzy languages:  these include  Arnoux-Rauzy languages on more than three letters, for which we give a (non-optimal) bound on the possible length of a clustering word, and {\it episturmian} languages, which include Sturmian languages, Arnoux-Rauzy languages, some periodic languages, and some intermediate cases which behave essentially like Sturmian words. Among episturmian languages on three letters we give in Theorem \ref{thepi} a full characterization of those which produce only finitely many clustering words: rather unexpectedly, these include not only Arnoux-Rauzy languages, but also some (not all) of the periodic and intermediate cases. 
 
 \section{Usual definitions}

 Let $\mathcal A$ be a finite set called the {\em  alphabet}, its elements being {\em letters}.
 A {\em word} $w$ of {\em length} $n=|w|$ is $a_1a_2\cdots a_n$, with $a_i \in {\mathcal A}$. The {\em concatenation} of two words $w$ and $w'$ is denoted by $ww'$.\\
 A word is  {\em primitive} if it is not a power of another word.\\
 The {\em reverse} of a word $w=w_1...w_n$
 is the word $\bar w=w_n...w_1$.\\
 
By a   language $\Lambda$ over $\mathcal A$ we mean  a {\it factorial extendable language}:
a collection of sets $(\Lambda_n)_{n\geq 0}$ where the only element of $\Lambda_0$ is the {\em empty word}, and where each $\Lambda_n$ for $n\geq 1$ consists of words of length $n$,  such that
for each $v\in \Lambda_n$ there exists $a,b\in \mathcal A$ with $av,vb\in \Lambda_{n+1}$,
and
each $v\in \Lambda_{n+1}$ can be written in the form $v=au=u'b$ with $a,b\in \mathcal A$ and $u,u'\in
\Lambda_n.$ \\
A word $v=v_1...v_r$  {\em occurs} at index $i$ in a word $w=w_1...w_s$  if $v_1=w_i$, ...$v_r=w_{i+r-1}$, we say also that $w$ contains $v$ and $v$ is a {\em factor} of $w$.\\

The {\em complexity} function of a language $\Lambda$ is $p(n)=\#\Lambda_n$, $n\geq 0$.\\
The {\em Rauzy graph} of length $n$ of a language $\Lambda$ is a directed graph whose vertex set consists of all words of length $n$ of $\Lambda$, with an edge from $w$ to $w'$ whenever $w=av$, $w'=vb$ for letters $a$ and $b$, and the word $avb$ is in $\Lambda$; this edge is then labelled by $b$.\\

 A word $w$ in $\Lambda$ is {\em  right special} (resp. {\em  left special}) if it has more than one {\em right extension} $wx$ (resp. {\em left extension} $xw$) in $\Lambda,$ with $x$ in $\mathcal A.$  If $w$ is both right special and
left special, then $w$ is {\em  bispecial}. If $\# \Lambda_1>1$, the empty word $\varepsilon$
 is bispecial.
 To {\em resolve} a bispecial word $w$ is to find all words in $\Lambda$ of the form $xwy$ for letters $x$ and $y$.\\
A {\em singular word} is $w=xvy$ for letters $x$, $y$,  such that some $x'vy$, $x'\neq x$, and $xvy'$, $y'\neq y$,  exist in $\Lambda$.\\

For a word $w$, we denote by $w^{\omega}$ the one-sided infinite word $www...$, and by $\Lambda_w$ the language consisting of all the factors of  $w^{\omega}$.
A language $\Lambda$ is {\em closed under reversal} if $w\in \Lambda \Leftrightarrow \bar w \in \Lambda$.\\
A language $\Lambda$ is {\em uniformly recurrent} if for every word $w$ in $\Lambda$, there exists a constant $K$ such that $w$ occurs in every word in $\Lambda$ of length at least $K$.

\section{Burrows-Wheeler and clustering}\label{sbw} 
Let  $\mathcal A=\{a_1<a_2< \cdots <a_r\}$ be an ordered alphabet. 

\begin{definition} The {\em (cyclic) conjugates} of $w$ are  the words $w_i\cdots w_nw_1\cdots w_{i-1}$, $1\leq i\leq n.$ If $w$ is primitive,  $w$ has precisely $n$ conjugates.  Let $w_{i,1}\cdots w_{i,n}$ denote the $i$-th  conjugate of $w$ where the $n$ conjugates of $w$ are ordered by ascending lexicographical order.\\
 Then the {\em  Burrows-Wheeler transform} of $w$, defined in \cite{bw} and denoted by $B(w)$, is the word $w_{1,n}w_{2,n}\cdots w_{n,n}.$
 In other words, $B(w)$ is obtained from $w$ by first ordering its conjugates in ascending order in a rectangular array, and then reading off the last column.\\
 We say $w$ is {\em clustering for the permutation $\pi$} if $B(w)=(\pi a_1)^{n_{\pi a_1}}\cdots (\pi a_r)^{n_{\pi a_r}}$, where $\pi\neq Id$ is a permutation on $\mathcal A$ and $n_a$ is the number of occurrences of $a$ in $w$ (we allow some of the $n_a$ to be $0$, thus, given the  order and $w$, there may be several possible $\pi$). We say $w$ is {\em perfectly} clustering if it is clustering for the {\em symmetric} permutation $\pi a_i=a_{r+1-i}$, $1\leq i\leq r$. \end{definition}

{\bf Non-primitive words}. As remarked in \cite{mrs}, the Burrows-Wheeler transform can be extended to a non-primitive word $w_1\cdots w_n$, by ordering its $n$ (non necessarily distinct) cyclic conjugates  by non-strictly increasing lexicographical order and taking the word made by their last letters. Then $B(v^m)$ is deduced from $B(v)$ by replacing each of its letters $x_i$ by $x_i^m$, and $v^m$ is clustering for  $\pi$ iff $v$ is clustering for  $\pi$.\\

We shall now relate clustering to an {\it order condition}. This condition can be traced to \cite{fz3}, but was first mentioned explicitly in \cite{zdl} in the particular case of symmetric permutations, and \cite{fhuz} in the general case, where it is studied extensively. 

 \begin{theorem}\label{car} For a given order $<$ on the alphabet $\mathcal A,$ a primitive word $w$ over $\mathcal A$ is clustering for the permutation $\pi$ if and only if every bispecial word
$v$ in the language $\Lambda_w$  satisfies the following {\em order condition}: whenever  $xvy$ and $x'vy'$ are in $\Lambda_w$ with letters  $x\neq x'$ and $y\neq y'$, then $\pi^{-1} x< \pi^{-1}  x'$ if and only if $y<y'$.  \\
Any bispecial word in $\Lambda_w$ is a factor of $ww$ and is of length at most $|w|-2.$ \end{theorem}
{\bf Proof}\\
We begin by proving the last assertion. Suppose $v$ is a bispecial of $\Lambda_w$. Then $v$ must occur at two different positions in some word $w^k$. If $|w|=n$ and $|v|\geq |w|-1$, this implies  in particular $w_i...w_nw_1...w_{i-2}=w_j...w_nw_1...w_{j-2}$ for $1<j-i<n$, and we notice that each $w_i$ is in at least one member of the equality, thus we get that $w$ is a power of a word whose length is the GCD of $n$ and $j-i$, which contradicts the primitivity. Thus the length of $v$ is at most $|w|-2$, and it occurs in $ww$.\\

We prove now that our order condition is equivalent to the following {\em modified order condition}: whenever $z=z_1...z_n$ and $z'=z'_1...z'_n$ are two different cyclic conjugates of $w$, $z<z'$ (lexicographically) if and only if $\pi^{-1}  z_k<\pi^{-1}  z'_k$ for the  largest $k\leq n$ such that $z_k\neq z'_k$.

Indeed, by definition $z<z'$ if and only if $z_j<z'_j$ for the  smallest $j\geq 1$ such that $z_j\neq z'_j$. If $w$ satisfies the order condition, we apply it  to the bispecial word $z_{k+1}...z_nz_1...z_{j-1}$, with $k$ and $j$ as defined, and get the modified order condition.

 Let $v$ be  a bispecial word in $\Lambda_w$; by the first paragraph of this proof
it can be written as $z_1...z_{k-1}$ for some $1\leq k\leq n$, with the convention that $k=1$ whenever $v$ is empty,  and at least two different cyclic conjugates  $z$ of $w$, and its possible extensions are the corresponding $z_nz_1...z_k$, thus, if  the modified order condition is satisfied, $v$  does satisfy the requirement of the order condition. \\

The modified order condition implies clustering, as then if two cyclic conjugates of $w$ satisfy $z<z'$,  their last letters $z_n$ and $z'_n$ satisfy either $z_n=z'_n$ or
$\pi^{-1} z_n<\pi^{-1} z'_n$.  

Suppose $w=w_1\cdots w_n$ is clustering for $\pi$. Suppose two cyclic conjugates of $w$ are such that $z_k\neq z'_k$, $z_j=z'_j$ for $k+1 \leq j\leq n$. Then $z<z'$ is (by definition of the lexicographical order) equivalent to $z_{k+1}...z_nz_1..z_k <z'_{k+1}...z'_nz'_1..z'_k$, and, as these two words have different last letters, because of the clustering this  is equivalent to $\pi^{-1} z_k<\pi^{-1}  z'_k$, thus the modified order condition is satisfied.
\qed\\

Theorem \ref{car} remains valid if $w=v^m$ is non-primitive (it can be slightly improved as there are less bispecial words to be considered, it is enough to look at factors  of $vv$ of length at most $|v|-2$).\\

The following consequences of Theorem \ref{car} or of \cite{fz4} seem to be new.

\begin{proposition} If $w$ clusters  for the order $<$ and the permutation $\pi$, its reverse clusters  for the $\pi$-order,  defined  by $x<_{\pi} y$ whenever
$\pi^{-1} x< \pi^{-1}  y$, and the permutation $\pi^{-1}$.
\end{proposition}
{\bf Proof}\\ 
This follows immediately from Theorem \ref{car}. \qed\\

\begin{proposition}\label{rev} Let $w$ be a word on  $\mathcal A$, ordered by $<.$.\\
If $w$ is perfectly clustering,  $\Lambda_w$ is closed under reversal. \\
If $\Lambda_w$ is closed under reversal, the following conditions are equivalent
\begin{enumerate}
\item $w$ is clustering.
\item $w$ is perfectly clustering.
\item For all words $u$ and $v$ with $u\neq \bar u$ and $v\neq \bar v,$  if $uv$ is conjugate to $w,$ then $u<\bar u$ if and only if $v<\bar v.$
\end{enumerate} 
\end{proposition}

{\bf Proof}\\
By Theorem 4 of \cite{fz4}, every perfectly clustering word $w$  is such that $ww$ is in the language $\Lambda$ generated by a minimal discrete interval exchange with the symmetric permutation (we refer the reader to \cite{fz4} for the definitions), and $\Lambda_w=\Lambda$.  It is known from \cite{fz3} that such a $\Lambda$  is stable under reversal, thus we get our first assertion. This could also be deduced from Corollary 4.4 of \cite{sp}. \\

We begin by showing the equivalence between $(1)$ and $(2)$ then we show that $(2) \Leftrightarrow (3).$  Clearly $(2)\Rightarrow (1).$
To see that $(1)\Rightarrow (2),$ assume that $w$ is clustering for some permutation $\pi$ on $\mathcal A.$ Let $\mathcal A'$ be the set of all letters $a\in \mathcal A$ which occur in $w.$ To show that $w$ is perfectly clustering, it suffices to show that $\pi^{-1} a < \pi^{-1}b \Leftrightarrow b<a$ for each pair of distinct letters $a,b \in \mathcal A'.$ To this end, we will show that the following set
\[\mathcal E=\{(a,b)\in \mathcal A' \times \mathcal A'\,: \, a\neq b\,\,\mbox{and}\,\, a<b \Leftrightarrow \pi^{-1}a <\pi^{-1}b \}\]
is empty.  We begin by establishing two claims:

{\bf Claim 1 :} Assume $xvy, x'vy' \in \Lambda_w$ with $v$ a word, letters  $x\neq x'$ and $y\neq y'.$ Then $(x,x')\in \mathcal E$ if and only if $(y,y') \in \mathcal E.$ 

{\bf Proof :} As $\Lambda_w$ is closed under reversal, we also have $y\bar v x, y'\bar v x' \in \Lambda_w.$ By Theorem~\ref{car} we obtain $y<y' \Leftrightarrow \pi^{-1}x < \pi^{-1} x'  \Leftrightarrow x<x'  \Leftrightarrow \pi^{-1}y < \pi^{-1}y'.$

{\bf Claim 2 :} Assume $xvy, x'vy' \in \Lambda_w$ with $v$  a word, letters $x\neq x'$ and $y\neq y'.$ If $(x,x')\in \mathcal E,$ then $x<x'\Leftrightarrow y<y'.$

{\bf Proof :} Again by Theorem~\ref{car} we have $x<x' \Leftrightarrow \pi^{-1}x <  \pi^{-1} x' \Leftrightarrow y < y'.$ \\

Now assume to the contrary that $\mathcal E \neq \emptyset $ and let $(x_1,y_1)\in \mathcal E.$ Without loss of generality we may assume that $x_1<y_1.$ Let $u$ and $v$ be conjugates of $w$ with $u$ beginning in $x_1$ and $v$ beginning in $y_1.$ Then we may write $u=x_1v_1x_2v_2\cdots x_nv_n$ and $v=y_1v_1y_2v_2\cdots y_nv_n$ for some $n\geq 2$ with words $v_i$ and letters $x_i\neq y_i$ for each $i=1,2,\ldots ,n.$ By application of Claims $1$ and $2$ we have that $(x_i,y_i)\in \mathcal E$ and $x_i<y_i$ for each one of the three rules is used infinitely many times, $i=1,2,\ldots ,n.$ As $u$ and $v$ are conjugate to one another,  in particular they have the same number of occurrences of each letter, and hence the same is true of the words $x=x_1x_2\cdots x_n$ and $y=y_1y_2\cdots y_n.$ Pick a permutation $\sigma$ of $\{1,2,\ldots ,n\}$ such that $y_i=x_{\sigma (i)}$ for each $i=1,2,\ldots ,n.$ It follows that $x_i<x_{\sigma (i)}$ for each $i=1,2,\ldots ,n.$ Putting $i$ equal to $\sigma^j(1)$ we obtain $x_{\sigma^j (1)}<x_{\sigma^{j+1}(1)}$ for each $j\geq 0.$ Thus $x_1<x_{\sigma(1)}<x_{\sigma^2(1)} < \cdots .$ Since $\sigma^{n!}(1)=1$ we  eventually get $x_1<x_1,$ a contradiction. 

We next show that $(2)\Rightarrow (3).$ So assume that $w$ is perfectly clustering. Then by Theorem~\ref{car} we have that $\Lambda_w$ satisfies the following order condition : whenever  $xzy$ and $x'zy'$ are in $\Lambda_w$ with  $x\neq x'$ and $y\neq y',$ we have $y<y' \Leftrightarrow  x'<x.$ Assume $uv$ is conjugate to $w$ with $u\neq \bar u$ and $v\neq \bar v.$ Write $u=rxtx'\bar r$ and $v=sy't'y\bar s$ with words $r, s, t, t'$, letters   $x,x',y,y',$ $x\neq x'$ and $y\neq y'.$ Thus $x'\bar r s y', y\bar s rx \in \Lambda_w$ and hence also $x\bar r s y \in \Lambda_w.$ Applying the order condition to the words $x\bar r s y$ and $x'\bar r s y'$ we obtain $y'<y \Leftrightarrow x<x'$ or equivalently $u<\bar u \Leftrightarrow v<\bar v$ as required.

Finally we show that $(3)\Rightarrow (2).$ Again by application of Theorem~\ref{car} it suffices to show that $\Lambda_w$ satisfies the following order condition : whenever  $xzy$ and $x'zy'$ are in $\Lambda_w$ with $|xzy|\leq |w|$,$z$ a word, $x,x', y, y'\in \mathcal A,$   $x\neq x'$ and $y\neq y',$ we have $y<y' \Leftrightarrow  x'<x.$
So assume that $xzy, x'zy' \in \Lambda_w$ with $|xzy|\leq |w|,$  $x\neq x'$ and $y\neq y'.$ Then $y'\bar z x' \in \Lambda_w$ since $\Lambda_w$ is closed under reversal. Let $w'$ be a conjugate of $w$ beginning in $xzy.$ If the words $xzy$ and $y'\bar z x'$  are equal, then in particular $x=y’$ and $y=x’,$ and therefore $y<y’$ iff $x’<x$, which is what we want. If they are different, then since $xzy$ and $y'\bar z x'$ cannot overlap one another,  we can write
$w'=xzyry'\bar z x' s$ for some choice of words $r,s.$ Put $u=zyry'\bar z$ and $v=x'sx.$ Then as $uv$ is conjugate to $w$ and $u\neq \bar u$ and $v\neq \bar v,$ we deduce $y<y' \Leftrightarrow u<\bar u \Leftrightarrow  v<\bar v \Leftrightarrow x'<x$ as required. \qed\\ 

We use now Theorem \ref{car} to give a simple criterion which will be useful to avoid clustering.

\begin{lemma}\label{ord} Let $x,y,z$ be three different letters in an alphabet $\mathcal 	A$, and  $w$ be a word on $\mathcal A$. Suppose $w$ is clustering for the order $<$ and the permutation $\pi$. Let $v$ be a bispecial word in $\Lambda_w$:
\begin{itemize} \item if  the  four words  $xvy$, $xvz$, $yvx$, $zvx$  are in $\Lambda_w$, then $x$ is not between $y$ and $z$ (or $z$ and $y$) for the order $<$, $x$ is not between $y$ and $z$ (or $z$ and $y$) for the order $<_{\pi}$, and $x$ is not on the same side of $y$ and $z$ for the orders $<$ and $<_{\pi}$ ;
 \item if  three of the  four words  $xvy$, $xvz$, $yvx$, $zvx$   are in $\Lambda_w$, then $x$ is not between $y$ and $z$ (or $z$ and $y$)  for the order $<$, or $x$ is not between $y$ and $z$ (or $z$ and $y$) for the order $<_{\pi}$. \end{itemize}\end{lemma}
{\bf Proof}\\
By Theorem \ref{car}, we have to check the order condition, for any fixed $<$ and $\pi$. To check the requirement of the order condition for the bispecial word $v$, we write the {\it extension graph} of $v$, with $x$ $y$ $z$ in the order $<$ on a line, $x$ $y$ $z$ in the order $<_{\pi}$ on a line below, and  an edge from $x'$ below to $y'$ above whenever $x'vy'$ is in $\Lambda_w$. If two of these edges have an intersection not reduced to an endpoint, the order condition is not satisfied. 

 In the first case, suppose $y<x<z$. if the $\pi$ order is $x,y,z$, or $x,z,y$, or $y,x,z$ the edges $xz$ and $zx$ intersect; if it is $z,x,y$, $xy$ and $zx$ intersect. All other possible $\pi$-orders are deduced by left/right symmetry, and no clustering is possible. The same is true if $z<y<x$ by left/right symmetry. 
 If $x<y$ and $x<_{\pi}y$, or  $y<x$ and $y<_{\pi}x$, $xy$ and $yx$ intersect whatever the position of $z$. Finally, the case where $x$ is the middle for the $\pi$-order is dealt with by up/down symmetry. 

Suppose now for example $xwy$, $ywx$ and $xwz$ are in $vv$, we have to test all orders where $x$ is twice in the middle. As above, if 
$x<y$ and $x<_{\pi}y$, or  $y<x$ and $y<_{\pi}x$, $xy$ and $yx$ intersect. There remain $z<x<y$ and $y<_{\pi}x<_{\pi}z$, and $y<x<z$ and $z<_{\pi}x<_{\pi}y$, and in both cases $xz$ intersects $yx$. And similarly for other sets of three words. 
\qed\\

\section{Arnoux-Rauzy}\label{sar}
\subsection{Definitions}\label{sdar}
Throughout  Section \ref{sar}, we use the alphabet $\{a,b,c\}$, which can be equipped  with any  one of the six possible orders. 

 \begin{definition}\label{dar} An {\em AR language} is a language on $\{a,b,c\}$ generated by
three families of words $A_{k}$, $B_{k}$, $C_{k}$, build recursively from $A_{0}=a$,
$B_{0}=b$, $C_{0}=c$, by using a
sequence of
combinatorial rules $(a)$, $(b)$, $(c)$, such that each one of the three rules is used infinitely many times, where
\begin{itemize}\item by rule $(a)$ at stage $k$, $A_{k+1}=A_{k}$,
$B_{k+1}=B_kA_k$, $C_{k+1}=C_kA_k$;
\item by rule $(b)$ at stage $k$, $A_{k+1}=A_kB_{k}$,
$B_{k+1}=B_{k}$, $C_{k+1}=C_kB_{k}$; \item
by rule $(c)$ at stage $k$, $A_{k+1}=A_kC_{k}$,
$B_{k+1}=B_kC_{k}$, $C_{k+1}=C_{k}$.
\end{itemize}
By an {\em AR word}  we shall mean a factor of an AR language. 

A {\em standard} AR word is an  $A_k$, $B_k$, or $C_k$, in an AR language
 
If the  rules at stage $k$ is $(x_k)$, $k\geq 0$, the word $D=x_0x_1....$ is called the {\em directive word} of $\Lambda$. The {\em Tribonacci language} is the AR language defined by the directive word $(abc)^{\omega}$.\end{definition}

Every AR language is uniformly recurrent and closed under reversal, and has one right special and one left special of each length \cite{ar}, thus AR languages are in the slightly more general class of episturmian languages, see Section \ref{sepi} below.

An equivalent way to define an AR language is through {\it AR morphisms}. For $x$, $y$ in $\{a,b,c\}$, we define $\sigma_xx=x$, $\sigma_xy=yx$ if $x\neq y$. For a word $w$, $\sigma_xw$ is defined  by making $\sigma_x$ a morphism for the concatenation, and the morphism $\sigma_w$ is defined to be $\sigma_{w_1}\circ \cdots \circ \sigma_{w_n}$ if $w=w_1...w_n$. We do the same with the morphisms $\tau_x$ defined by $\tau_xy=\overline {\sigma_xy}$ for each $x,y\in \mathcal A.$
If the directive word of $\Lambda$ is $D=x_0x_1...$, we put $D_k=x_0...x_{k-1}$, and we have $A_k=\sigma_{D_k}a$, 
 $B_k=\sigma_{D_k}b$,  $C_k=\sigma_{D_k}c$, $\bar A_k=\tau_{D_k}a$, 
 $\bar B_k=\tau_{D_k}b$,  $\bar C_k=\tau_{D_k}c$. Being closed under reversal, $\Lambda$ can be generated either by $\sigma_{D_k}x$, $x \in \mathcal A$, $k\geq 0$, or by  $\tau_{D_k}x$, $x \in \mathcal A$, $k\geq 0$.  \\

For an AR language $\Lambda$, note first that for all $k$ $A_k$ begins with $a$, $B_k$ with $b$, $C_k$ with $c$. As explained in \cite{ar}, the three rules correspond to the building of the successive bispecials $w_k$ in $\Lambda$, with rules $(a)$, $(b)$, $(c)$ corresponding respectively to $w_{k+1}=w_kA_k$, $w_{k+1}=w_kB_k$, $w_{k+1}=w_kC_k$, starting with $w_0$  being the empty word. From this and the closure under reversal, we deduce that AR rule $(j)$,  $j=a,b,c$,  is used at stage $k$ if and only if 
the bispecial $w_k$ is resolved by  $\{aw_kj,
bw_kj, cw_kj, jw_ka, jw_kb, jw_kc\}$, two of these six words being equal. Moreover, the words $A_k$, $B_k$, $C_k$ are the {\it (suffix) return words} of $w_k$ i.e. $w_kZ_k$ contains $w_k$ as a prefix and suffix and at no other place, for $Z=A,B,C$. \\

An AR word $w$ belonging to an AR language whose first rule is $(x)$ will be such that each letter of $w$ which is not $x$ is preceded (except if it is the first letter of $w$), and followed (except if it is the last letter of $w$) by $x$, and $x$ is the only letter with this property. We call $x$ the {\it separating letter} of $w$. 

 We recall the  description of the Rauzy graphs for AR languages from \cite{ar}: there are a left special factor $G$ with three left extensions, a right special factor $R$ with three right  extensions, a central branch (with at least  one vertex) from $G$ to $R$, and three branches from $R$ to $G$. The three {\it elementary circuits} in the  Rauzy graphs of length $n$ begin at 
 $R$ and follow one of the three branches from $R$ to $G$ then the central branch. Their labels are $A_k$, $B_k$, $C_k$ for $|w_{k-1}|+1\leq n \leq |w_k|$.\\

We shall always use the obvious notation that if $x$, $y$ or $z$  is $a$, $b$ or $c$, $X$, $Y$  and $Z$ are the corresponding $A$, $B$ or $C$.\\

Let $(Rabc)$ be the following assumption: {\it the rule at stage $0$ is $(a)$, and the first rule different from $(a)$ is $(b)$}. If it is not satisfied, we can make a permutation on the letters. When $(Rabc)$ holds,  we define $\lambda_1>0$ as the stage of the first rule $(b)$, $\lambda_2>\lambda_1$ as the stage of the first rule $(c)$.\\

The following {\it LMS notation} is defined in \cite{cc}; it is equivalent to the ABC notation, and will be useful to express and show some of our results.  

\begin{lemma}\label{lms} Assuming $(Rabc)$,
for all $k>\lambda_1$, $A_k$, $B_k$, and $C_k$ have three different lengths, and we rename them such that $|S_k|<|M_k|<|L_k|$. We put  $S_k=A_k$, $M_k=C_k$, $L_k=B_k$ for all $1\leq k\leq \lambda_1$, $S_0=c$, $M_0=b$, $L_0=a$. The AR rules can be written as 
\begin{itemize}\item if $w_{k+1}=w_kS_k$, $S_{k+1}=S_{k}$,
$M_{k+1}=M_kS_k$, $L_{k+1}=M_kL_k$;
\item if $w_{k+1}=w_kM_k$,, $S_{k+1}=M_k$,
$M_{k+1}=S_{k}M_k$, $L_{k+1}=L_kM_{k}$; \item
if $w_{k+1}=w_kL_k$, $S_{k+1}=L_k$,
$M_{k+1}=S_{k}L_k$, $L_{k+1}=M_kL_{k}$.
\end{itemize}
 
We have $w_{p+1}=w_pL_p$ whenever $p=0$, $p=\lambda_1$, $p=\lambda_2$, or at stage $p$ we have  a rule $(x_1)$ preceded by a string of rules $(x_2)$ and a string of rules $(x_3)$, for $\{x_1,x_2,x_3\}=\{a,b,c\}$, thus this happens for infinitely many $p$.
 We have  $w_{p+1}=w_pS_p$ whenever the rules at stages $p-1$ and $p$ are the same, and  $w_{p+1}=w_pM_p$ for the remaining $p$. 
 \end{lemma}
{\bf Proof}\\ We have $w_k=a^k$,  $A_k=a$, $B_k=ba^k$, $C_k=ca^k$ for all $1\leq k\leq \lambda_1$, then
$w_{\lambda_1+1}=a^{\lambda_1}ba^{\lambda_1}$, and
 $A_{\lambda_1+1}=aba^{\lambda_1}=M_{\lambda_1+1}$,  $B_{\lambda_1+1}=ba^{\lambda_1}=S_{\lambda_1+1}$
 $C_{\lambda_1+1}=ca^{\lambda_1}ba^{\lambda_1}=L_{\lambda_1+1}$ have three different lengths; this is preserved by further rules. The other assertions are straightforward, see \cite{cc} for more details. \qed\\

\subsection{Lengths, squares, conjugates}
\begin{lemma}\label{rel} Assuming $(Rabc)$,\\
$A_p$ is a suffix of $w_p$  iff $|A_p|\leq |w_p|$ iff $p\geq 1$,\\
$B_p$ is a suffix of $w_p$ iff $|B_p|\leq |w_p|$  iff $p\geq \lambda_1+1$,\\
$C_p$ is a suffix of $w_p$ iff $|C_p|\leq |w_p|$  iff $p\geq \lambda_2+1$,\\
 $|M_p|+|S_p|>|L_p|$ for all $p>0$,\\
 $|L_p|<|M_{p+1}|$ for all $p>0$.
\end{lemma}
{\bf Proof}\\
By the analysis of Lemma \ref{lms}, $A_p$ is a suffix of $w_p$ for $p=1$ and strictly longer than $w_p$ for $p=0$,  $B_p$ is a suffix of $w_p$ for $p=\lambda_1+1$ and strictly longer than $B_p$ for $p\leq \lambda_1$. $C_p$ is strictly longer than   $w_p$ for $p= \lambda_1+1$. If $C_p$ is strictly longer than $w_p$ for some $p\leq \lambda_2-1$, then $w_{p+1}=w_pY_p$ and $C_{p+1}=C_pY_p$ for $Y=A,B$, thus  $C_{p+1}$ is strictly longer than $w_{p+1}$. 
Then $w_{\lambda_2+1}=w_{\lambda_2}C_{\lambda_2}$ and $C_{\lambda_2+1}=C_{\lambda_2}$ is a suffix of $w_{\lambda_2+1}$. Suppose now $Z_p$ is a suffix of $w_p$, for $Z=A,B,C$. Then either $w_{p+1}=w_pZ_p$ and $Z_{p+1}=Z_p$, or 
$w_{p+1}=w_pY_p$ and $Z_{p+1}=Z_pY_p$, thus in both cases $Z_{p+1}$ is a suffix of $w_{p+1}$, and thus shorter. 

 The fourth  assertion  is proved in \cite{cfme}, although with non-strict inequalities, but the proof, using the AR rules, does give the strict ones.

 The last one comes from the fourth one and the fact that $|M_{p+1}|$ is either $|M_p|+|S_p|$ or  $|M_p|+|L_p|$.\qed\\

\begin{proposition}\label{sq} 
In an AR language $\Lambda$ satisfying $(Rabc)$, the primitive words $v$  such that $vv$ is in $\Lambda$ are all the $A_p$ for $p\geq 0$,  the $B_p$ for $p\geq \lambda_1$, the $C_p$ for $p\geq \lambda_2$, and some (possibly none) of their (cyclic) conjugates. \\
For a primitive word $v$, the following assertions are equivalent 
\begin{itemize}
\item $v$ is conjugate to a standard AR word,
\item $vv$ is an AR word,
\item  $v'v'$ is an AR word for some conjugate $v'$ of $v$,
\item   all the conjugates of $v$ are in some AR language $\Lambda$,
\item $x$ can be de-substituted down (using the six AR morphisms $\tau_x$ and $\sigma_x$,  $x \in \mathcal A$) to a single letter. 
 \end{itemize}\end{proposition}
{\bf Proof}\\ 
Let $Z_p=A_p$ for $p\geq 1$, $B_p$ for $p\geq \lambda_1+1$, or $C_p$ for $p\geq \lambda_2+1$. Then $w_pZ_p$ is in $\Lambda$ while by Lemma \ref{rel} $Z_p$ is a suffix of $w_p$, thus $Z_p^2$ is in $\Lambda$. This is true also for $A_0=a$
as $aa$ is a suffix of $w_1A_1$, for  $B_{\lambda_1}$ as 
$w_{\lambda_1+1}B_{\lambda_1+1}=w_{\lambda_1}B_{\lambda_1}^2$ is in $\Lambda$, and for $C_{\lambda_2}$ as 
$w_{\lambda_2+1}C_{\lambda_2+1}=w_{\lambda_2}C_{\lambda_2}^2$ is in $\Lambda$. Note that some conjugates of the $Z_p$ may have the same property, but for the Tribonacci language and $v=A_2=aba$, $vv$ is in $\Lambda$ but no $v'v'$ for $v'=baa$ or $v'=aab$.

As remarked in \cite{brl}, for a primitive $v$, if $vv$ is in $\Lambda$, then $v$ is the label of a circuit in the Rauzy graph of length $|v|$ where no vertex is used more than once; thus $v$ can only be some conjugate of some $A_p$, $B_p$, or $C_p$, for a $p$ such that
 $|w_{p-1}|+1\leq |v| \leq |w_p|$. This will not be the case for the conjugates of $B_p$, $p<\lambda_1$, as then $B_p$ is strictly longer than $w_p$ and $B_{p+1}\neq B_p$ because the rule at stage $p$ is not $(b)$, nor  for the conjugates of $C_p$, $p<\lambda_2$, as then $C_p$ is strictly longer than $w_p$ and $C_{p+1}\neq C_p$ because the rule at stage $p$ is not $(c)$. Thus we get the first assertion.\\

Suppose  a primitive $v$ is conjugate to a standard AR word; to show that $vv$ is always an AR word, what remains to prove is that  when   $v$ is one of the initial $A_p$, $B_p$, $C_p$, or is conjugate to any $A_p$, $B_p$, $C_p$ in an AR language $\Lambda$ satisfying $(Rabc)$, then $vv$ is in an AR language $\Lambda'$. This is true for $B_0$, $C_0$ by exchanging $a$ with $b$ or $c$. For $B_p$, $1\leq p\leq \lambda_1$, this will be true for $\Lambda'$ defined by the same rules as $\Lambda$ up to stage $p-1$, then rule $(b)$ at stage $p$, and any admissible sequence  of rules beyond: in this language our $B_p$ is $B'_{\lambda'_1}$. For 
$C_p$, $1\leq p\leq \lambda_1$, this will be true for $\Lambda"$ deduced from the previous  $\Lambda'$ by exchanging $b$ and $c$.
 For $C_p$, $\lambda_1+1\leq p\leq \lambda_2-1$, this will be true for $\Lambda'$ defined by the same rules as $\Lambda$ up to stage $p-1$, then rule $(c)$ at stage $p$, and any admissible sequence  of rule beyond: in this language our $C_p$ is $C'_{\lambda'_1}$. 
Let now $u$ be a conjugate to some $Z_p$ which is a suffix of $w_p$; then we define $\Lambda'$  by the same rules as $\Lambda$ up to stage $p-1$, then rule $(z)$ at stage $p$, and any admissible sequence  of rules beyond: then $w_pZ_p^2$, thus $Z_p^3$, thus $uu$ is in $\Lambda$. This will still be true for $Z_p=B_{\lambda_1}$ or $Z_p=C_{\lambda_2}$  if we define $\Lambda'$ by  the same rules as $\Lambda$ up to stage $p$, then rule $(z)$ at stage $p+1$, and any admissible sequence  of rules beyond.

If $vv$ is in $\Lambda$, then $v'v'$ is in $\Lambda$ for the conjugate $v$, and if $v'v'$ is in $\Lambda$ all the  conjugates of $v'$, thus of $v$,  are in $\Lambda$. 

Suppose all the  conjugates of $v$ are in an AR language $\Lambda$. In the Rauzy graph  of length $|v|$ of $\Lambda$, we see each  conjugate $v^{(i)}$; there is an edge between $v^{(i)}$ and the next one in the circular order, because either $v^{(i)}$ has only one right extension or $v^{(i)}$ has all the possible right extensions.  Thus there is a circular path  whose vertices are all the $v^{(i)}$, each one occurring only once as $v$ is primitive; at least one vertex $v'$ in this path is on the central branch, and if we start from this point the circular path is an allowed path. Thus $v'v'$ is in $\Lambda$, as remarked in \cite{brl}, thus  $v'$ is conjugate to a standard AR word by the reasoning above, thus $v$ is conjugate to a standard AR word. Thus we have proved the equivalence of the first three items in the second assertion. \\

To deal with the last item, let us show first that a standard AR word $Z$  is conjugate to its reverse. This is  true if $Z$ has one letter. Other $Z$ are of the form $Z=\sigma_xZ'$, for a shorter standard AR word $Z'$. Then, if we suppose $Z'$ is conjugate to $\bar{Z'}$, we get that $Z=\sigma_xZ'$ is conjugate to $\sigma_x\bar{Z'}$, and the latter is conjugate to $\bar Z=\tau_x\bar{Z'}$ as these  two words are of the form $yw$ and $wy$ for some letter $y$. 

Let $v$ be a standard AR word, with separating letter  $a$, and  $v'$  a conjugate of $v$. Then $v'$ must either begin or end in $a$, otherwise $v’ =xv”y$ with each of $x$ and $y$ different from $a$; then every conjugate of $v'$ which begins in $a$ must contain $yx$ as a factor, and in particular $v$ contains $yx$ which contradicts the fact that  $a$ is the separating letter of $w$.  If $v'$ ends in $a$, then $v'=\tau_au'$, by the properties of the separating letter, and  if $v'$  end in $a$ then  $v'=\sigma_au'$. 
If $v'=\sigma_u'$, then $u'$ is conjugate to the standard AR  word $u$ such  that $v=\sigma_au$. If $v'=\tau_au'$, then $v'$ is also conjugate to $\bar v$, $\bar v=\tau_a \bar u$ $u'$ for a  standard AR  word $u$, and  $u'$ is conjugate to $u$. Then we apply the same process to $u'$ as long as $u'$ has at least two letters, and end when w get to a single letter.

Finally, to get that the last item implies the first, it is enough to prove that if $u'$ is conjugate to a standard AR word $u$, then, for any letter $x$, $\sigma_xu'$ or $\tau_xu'$ is also conjugate to a standard AR word, and this is true as both these words are conjugate to $\sigma_xu$. 
 \qed\\

For non-primitive words Proposition \ref{sq} fails: in  an AR language $\Lambda$ where there are five consecutive rules $(a)$ at stages $p$ to $p+4$, $A_p^4$ is in $\Lambda$ while $A_p^2$ is not any $A_k$, $B_k$ or $C_k$ in $\Lambda$. If we take  $p=0$, we see that $A_0^2=aa$ is not an $A_k$, $B_k$ or $C_k$ in any AR language. 

For $v$ to be conjugate to $v'$ such that $v'v'$ is an  AR word, it is not enough that all the conjugates of $v$  are AR words, take $abc$ for example. Also, a conjugate of a standard AR word is not necessarily a standard AR word, for example $caba$ is standard in the Tribonacci language, but $abac$ is not standard in any AR language, otherwise it could be written as the concatenation of two words with the same last letter.

  To emphasize the part played by the cyclic conjugates, we notice the following fact.

\begin{proposition} Let $\Lambda$ be any language on an alphabet $\mathcal A$ such that for every $w$ in $\Lambda$ and $a$ in $\mathcal A$, there exists $v$ in $\Lambda$ such that $wva$ is in $\Lambda$. Then the closure of $\Lambda$ for the cyclic conjugacy is made of all the possible words on $\mathcal A$.\end{proposition}
{\bf Proof}\\ Let $x_1x_2…x_n$ be any word. We will show $x_1x_2…x_n$  is in $\Lambda'$.
Let $x_nwx_{n-1}$ be in $\Lambda$  such that $w$ contains  each letter sufficiently many times.
So $x_{n-1}x_nw$ is in $\Lambda'$. Now write $w=ux_{n-2}v$ where $u$  contains each letter sufficiently many times. So $x_{n-1}x_nux_{n-2}$ is in $\Lambda'$ (factorial property of languages), hence $x_{n-2}x_{n-1}x_nu$ is in $\Lambda'$ (closed under cyclic conjugates). And so on for $x_{n-3}$, $x_{n-4}$, ... \qed\\

\subsection{Non-clustering AR words}
\begin{lemma}\label{caro}
Let $\Lambda$ be an AR language satisfying $(Rabc)$. If $u\in\Lambda$ contains  at least three of the  non-singular words $xw_{\lambda _2}y$, $x,y$ in  $ \{a,b,c\}$, $(x,y)\neq (c,c)$,  $u$ cannot cluster for any permutation $\pi$ and any order on $\{a,b,c\}$, nor can $v$ if $u=vv$.
 \end{lemma}
{\bf Proof}\\
Suppose $u$ contains at least three of $cw_{\lambda_2}a$,  $cw_{\lambda_2}b$ ,
$bw_{\lambda_2}c$,  $aw_{\lambda_2}c$. As $w_{\lambda_1}$ is both a prefix and a suffix of $w_{\lambda_2}$, $u$ contains
 also  $cw_{\lambda_1}$,  $w_{\lambda_1}c$  and at least three of  $aw_{\lambda_1}$,  $w_{\lambda_1}a$,  $bw_{\lambda_1}$,  $w_{\lambda_1}b$. As $u$ is in $\Lambda$ and the rule at stage $\lambda_1$ is $(b)$, $u$ contains
$cw_{\lambda_1}b$,  $bw_{\lambda_1}c$ and at least either $aw_{\lambda_1}b$ or  $bw_{\lambda_1}a$. As  $w_{\lambda_1}$ begins and ends with $a$, $u$ contains $ca$, $ac$, $ab$, $ba$, $aa$. 

Suppose $u$ clusters, then we apply Lemma \ref{ord}. By its first assertion applied to the empty bispecial, we must assign $a$ to an end of the $<$ order and the opposite one of $<_{\pi}$. Then, by the second assertion applied to $w_{\lambda_1}$ and $w_{\lambda_2}$, as we cannot give two ends to both $b$ and $c$, we must assign $b$ to the middle of one order, and $c$ to the middle of the  other one. Thus, up to left/right and up/down symmetries in the pictures in the proof of Lemma \ref{ord}, we have $a<b<c$ and $b<_{\pi}c<_{\pi}a$. But we know that $cw_{\lambda_1}b$ and $bw_{\lambda_1}c$ are factors of $u$, which gives two intersecting lines $bc$ and $cb$ in the picture for $w_{\lambda_1}$, and this contradicts Theorem \ref{car}.\qed\\

Together with uniform recurrence,  Lemma \ref{caro} provides our main answer to Dolce's question: {\it in a given AR language, there are only finitely many clustering words}. This is noticed in \cite{bal} in the particular case of an infinite sequence of words in the Tribonacci language called the  Tribonacci standard words, for which the Burrows-Wheeler transform is explicitly computed, and this is generalized to $r$-Bonacci, see Section \ref{slga} below.  Note also that both AR words and clustering are mentioned in \cite{sp}, but no relation between these notions is established.

The following theorem gives an estimation for the maximal length of a clustering word, with some claims to optimality. It relies on a method used in \cite{hm} for Sturmian languages and \cite{cc} for AR languages.

\begin{theorem}\label{arc}
We recall that, assuming $(Rabc)$, 
 $\lambda_2$ is  the stage of  the first rule $(c)$; Let $\lambda_a$  be the stage of the last rule $(a)$ before $\lambda_2$, $\lambda_b$  the stage of the last rule $(b)$ before $\lambda_2$,  $\mu_a$  the stage of the first rule $(a)$ after $\lambda_2$,  $\mu_b$  the stage of the first rule $(b)$ after $\lambda_2$. Let $x$ and $y$ be the elements of $\{a,b\}$ such that $\mu_x<\mu_y$.\\
Then no word of length at least $|w_{\lambda_y}|+\max(|C_{\mu_y+1}|, |X_{\mu_y+1}|)+1$ can cluster for any permutation  and any order.
The bound is sharp for the Tribonacci language.\\
The bounds  are not sharp in general, both in the case $\lambda_x<\lambda_y$ and in the case $\lambda_y<\lambda_x$. 
\end{theorem}
{\bf Proof}\\ 
We want to estimate the minimal length of a word containing at least three $xw_{\lambda _2}y$,  $(x,y)\neq (c,c)$, so that we can apply Lemma \ref{caro}.
As is noticed in \cite{hm}, a word $w$ occurs in any word of $\Lambda$ whose length is  at least $|w|-1+t(w)$, where $t(w)$ is the maximal return time of $w$, i.e. the  maximal possible difference between the indexes of two consecutive occurrences of $w$; and $t(w)$ is the same as $t(u)$, where $u$ is the longest singular word contained in $w$, or $u$ is a single letter if $w$ contains no singular word.

We reprove, with other notations, Lemma 2.3 of \cite{cc}. Let $v=zw_{p-1}z$ be a singular word. We define four assertions:
\begin{itemize}
\item[$(OA_q)$] $v$ occurs once in $w_qA_q$, $v$ does not occur in $w_qB_q$ or $w_qC_q$,
\item[$(OB_q)$] $v$ occurs once in $w_qB_q$, $v$ does not occur in $w_qA_q$ or $w_qC_q$,
\item[$(OC_q)$] $v$ occurs once in $w_qC_q$, $v$ does not occur in $w_qA_q$ or $w_qB_q$,
\item[$(OT_q)$] $v$ occurs at least once in $w_qA_q$, $w_qB_q$ and  $w_qC_q$, the maximal return time of $v$ is $|L_q|=\max(|A_q|,|B_q|,|C_q|)$.
\end{itemize}

If $w_p=w_{p-1}Z_{p-1}$,  $z$ is the first letter of $Z_{p-1}$ and $Z_p=Z_{p-1}$, thus  $(OZ_p)$ holds. Then the AR rules imply that
 if the  rule at stage $q$ is $(a)$, $(OA_q)$ 
implies $(OT_{q+1})$, $(OB_q)$ 
implies $(OB_{q+1})$, $(OC_q)$
implies $(OC_{q+1})$, and mutatis mutandis for rules $(b)$ and $(c)$.\\

We need to know the maximal return times of 
$u=yw_{\lambda_y}y$ and $u'=xw_{\lambda_x}x$. The analysis above implies that $t(u)$ is known as soon as  we see a rule $(y)$ after $\lambda_2$, which happens at stage $\mu_y$, and that $t(u)=|L_{\mu_y+1}|$, and similarly $t(u')=|L_{\mu_x+1}|$.\\

The rule at stage $\lambda_2$ is $(c)$. There are only rules $(c)$ (or none) (strictly) between $\lambda_2$ and $\mu_x$, there are only rules  $(x)$ or $(c)$ (or none) (strictly) between $\mu_x$ and $\mu_y$. By Lemma \ref{lms}, for $p=\lambda_2$,  and $p=\mu_y$, we have $w_{p+1}=w_pL_p$; this is true also for $p=\mu_x$ if $\lambda_x<\lambda_y$, and in any case this happens for no other $\lambda_2\leq p\leq \mu_y$. In particular, we get that $L_{\mu_y}=Y_{\mu_y}$ and $L_{\mu_y+1}$ is the longest 
of $X_{\mu_y+1}$ and $C_{\mu_y+1}$.\\  

{\bf First case,  $\lambda_x<\lambda_y$}. \\
Then $\lambda_y=\lambda_2-1$. We know that $xw_{\lambda_2}c$ and $cw_{\lambda_2}x$,  which contain $u'$ as maximal singular word,  occur in any word  in $\Lambda$ of length at least $|w_{\lambda _2}| + |L_{\mu_x+1}|+1$, and this is smaller than the required bound as $|w_{\lambda_2}|-|w_{\lambda_2-1}|=|Y_{\lambda_2-1}|<|L_{\lambda_2}|<|M_{\mu_y}|\leq |L_{\mu_y+1}|-|L_{\mu_x+1}|$.

We know also that $yw_{\lambda_2}c$ and $cw_{\lambda_2}y$,  which contain $u$ as maximal singular word,  occur in any word  in $\Lambda$ of length at least $|w_{\lambda _2}| + |L_{\mu_y+1}|+1$, but we can improve this bound a little if we want only to see $yw_{\lambda_2}c$ or $cw_{\lambda_2}y$.  Indeed, $u$  occurs in any word $Z$ in $\Lambda$ of length $|L_{\mu_y}|+|u|-1$. Also,  $u$ is a prefix of $yw_{\lambda_2}c$ and there is only one way to extend $u$ to the right to length $|w_{\lambda_2}|+2$, giving $yw_{\lambda_2}c$, $u$ is a suffix of $cw_{\lambda_2}y$, and  there is only one way to extend $u$ to the left to length $|w_{\lambda_2}|+2$, giving $cw_{\lambda_2}y$. Thus $yw_{\lambda_2}c$ or $cw_{\lambda_2}y$ is in $Z$, provided $Z$ is long enough to ensure that we can extend $u$ to the right or left as far as that length while remaining in $Z$; in the worst case, we can extend it by a length $\frac{|Z|-|u|}2$, so we have to check this is at least $|w_{\lambda_2}|+2-|u|$. Thus we have to prove that $|L_{\mu_y+1}|>2(|w_{\lambda_2}|-|w_{\lambda_2-1}|)$: the right side is   $2|Y_{\lambda_2-1}|$  while the left side is at least $|L_{\lambda_2+1}|=|X_{\lambda_2-1}|+|C_{\lambda_2-1}|+2|Y_{\lambda_2-1}|$. And we conclude by Lemma \ref{caro}.\\

{\bf Second case,  $\lambda_y<\lambda_x$}.\\ Then $\lambda_x=\lambda_2-1$. Again,  we have to check that 
$xw_{\lambda_2}c$ and $cw_{\lambda_2}x$ occur in any word of the required length, and $yw_{\lambda_2}c$ or $cw_{\lambda_2}y$ occur in any word of the required length. Using the same methods as  in the previous paragraph, this is done by  checking that 
$|w_{\lambda_2}|-|w_{\lambda_y}|<|L_{\mu_y+1}|-|L_{\mu_x+1}|$ and 
$|L_{\mu_x+1}|>2(|w_{\lambda_2}|-|w_{\lambda_y}|)$. We have  $|w_{\lambda_2}|-|w_{\lambda_y}|=t(|Y_{\lambda_y}|+|X_{\lambda_y}|)$, for some positive  integer $t$. Knowing  the rules between  stages $\lambda_2$ and  $\mu_x$, we get that  both $|L_{\mu_x+1}|$ and
$|M_{\mu_x+1}|$ are   at least $2t(|X_{\lambda_y}|+|Y_{\lambda_y}|)$. Then we can conclude, using also that $|L_{\mu_y+1}|-|L_{\mu_x+1}|$ is at least $|M_{\mu_y}|\geq |M_{\mu_x+1}|$.\\

Suppose the  directive word of $\Lambda$ begins with $abcab$, as for Tribonacci. Then $\lambda_a<\lambda_b$, $x=a$, $y=b$, and the bound above is $26$. The primitive  word $aB_4C_4=aL_4M_4=(abacaba)^2cabaabacaba$, of length $25$, is in $\Lambda$ and does cluster, perfectly, for the order $a<c<b$ or the order $b<c<a$: this can be checked by hand, but is also a consequence of Proposition \ref{long} below, see Corollary \ref{sha} which also gives other cases where the bound is sharp. A more detailed study of the Tribonacci case is in Example \ref{e1} below.\\

Suppose the directive word of $\Lambda$  begins with $abacba$, we have $\lambda_b<\lambda_a$, $x=b$, $y=a$.  The bound in the first assertion is $45$, and we look at words of length $44$. By the above  reasoning  $bw_{\lambda_2}c$ and $cw_{\lambda_2}b$ must occur in all words of this length; the assertion $OT_{\mu_a+1}$ above holds for the singular word $aw_{\lambda_2-1}a$, and, by looking precisely at its occurrences, we check that the only word  in $\Lambda$ of length $44$ without $aw_{\lambda_2}c$ and without $cw_{\lambda_2}a$ is $v=w_{\lambda_a}L_{\mu_a}M_{\mu_a}$. We check that $v$ cannot cluster for any order and permutation, by hand or see Example \ref{e2} below, nor can any word of length $44$.\\

Suppose the directive word of $\Lambda$ begins with $abcba$, we have $\lambda_a<\lambda_b$, $x=b$, $y=a$.  The bound is $24$, and by the same reasoning as in last paragraph the only word of length $23$ which might cluster is $v'=w_{\lambda_a}L_{\lambda_3}M_{\lambda_3}$.  We check that $v'$ cannot cluster for any order and permutation, by hand or see Example \ref{e3} below, nor can any word of length $23$. \qed\\

\subsection{AR words conjugate to standard}
\begin{lemma}\label{rec} If  a bispecial $v$  in a language $\Lambda_w$ is resolved by a subset of $\{avb,avc,ava,bva,cva\}$, $v$   satisfies the requirement of the order condition for any order $<$ such that $a$ is at an end, and the symmetric permutation.  If  in  $\Lambda_w$ a bispecial $v'$ is resolved by a subset of $\{av'b,cv'b,bv'b,bv'a,bv'c\}$, both $v$ and $v'$  satisfy the requirement of  the order condition for the  orders $a<c<b$ or $b<c<a$, and the symmetric permutation. \end{lemma}
{\bf Proof}\\
We draw the extension graphs as in the proof  of Lemma \ref{ord} and check that any two edges do not intersect except at their endpoints. \qed\\

Lemma \ref{rec} provides a partial converse to Lemma \ref{caro}, as it allows us to build clustering AR words in the absence of the obstructions in its hypothesis, but it does not give  a necessary condition for clustering, as we shall see in Section \ref{sns} below.

\begin{proposition}\label{list}  
With the assumption $(Rabc)$ and the  notations of Theorem \ref{arc} above, 
$Y_p$ clusters if and only if $p\leq \mu_x$,
$C_p$ and $X_p$ cluster if and only if  $p\leq \mu_y$.  
\end{proposition}
{\bf Proof}\\
We begin by the negative direction. Let $Z_p$ be $A_p$, $B_p$ or $C_p$ for the values in the hypotheses. Then, using the rules between $\lambda_2$ and  $\mu_y$ as determined in the proof of Theorem \ref{arc}, we check that $A_{\lambda_2}$, $B_{\lambda_2}$ and $C_{\lambda_2}$  all appear in the decomposition of $Z_p$ by the AR rules; as each $A_{\lambda_2}$ and $B_{\lambda_2}$ in this decomposition is preceded by $w_{\lambda_2}$, if $|Z_p|\geq |w_{\lambda_2}|$ then $w_{\lambda_2}A_{\lambda_2}$ and $w_{\lambda_2}B_{\lambda_2}$, thus $w_{\lambda_2}a$ and $w_{\lambda_2}b$, occur in $Z_p^2$; indeed, for these values of $p$, we have  $|Z_p|\geq |w_{\lambda_2}|+1$ by Lemma \ref{rel}, thus $cw_{\lambda_2}a$ and $cw_{\lambda_2}b$ occur in $Z_p^2$. A symmetric reasoning holds  for $aw_{\lambda_2}c$ and $bw_{\lambda_2}c$,  as $w_{\lambda_2}A_{\lambda_2}=A'_{\lambda_2}w_{\lambda_2}$, 
 $w_{\lambda_2}B_{\lambda_2}=B'_{\lambda_2}w_{\lambda_2}$, where $A'_{\lambda_2}$ ends with $a$ and $B'_{\lambda_2}$ ends  with $b$. This contradicts the clustering by Lemma \ref{caro}.\\

In the positive direction, let $Z_p$ be an $A_p$ for $p\geq 0$, or a $B_p$ for $p\geq \lambda_1$, or  a $C_p$ for $p\geq \lambda_2$.
By the reasoning of Proposition \ref{sq}, $Z_p^2$ is a suffix of $w_pZ_p$, or $w_{p+1}Z_{p+1}$ if $Z_p=A_0$,  $Z_p=B_{\lambda_1}$ or $Z_p=C_{\lambda_2}$, and
 the bispecials in the language $\Lambda_{Z_p}$ are resolved by AR rules. By Lemma \ref{rec}, those which are resolved by rule $(a)$ or $(b)$ satisfy the order condition ; as for bispecials $w_t$ resolved by rule $(c)$, they do satisfy the order condition if they are resolved in $\Lambda_{Z_p}$ by $\{cwc, awc, cwa\}$ or $\{cwc, bwc, cwb\}$. This will happen if the longest singular word $aw_{q'}a$ occurring in $awc$ and $cwa$, or  the longest singular word $bw_qb$ occurring in $cwb$ and $bwc$, does not occur in $Z_p^2$.

 Using the rules between $\lambda_2$ and  $\mu_y$ as determined  in the proof of Theorem \ref{arc}, we track  $bw_{\lambda_b}b$ and $aw_{\lambda_a}a$ as in Theorem \ref{arc}, and get that 
one of them, namely $u=yw_{\lambda_y}y$, does not occur in $w_pZ_p$, nor in $w_{p+1}Z_{p+1}$ when needed, hence in $Z_p^2$, and thus $w_{\lambda_2}$ satisfies the order condition in $\Lambda_{Z_p}$. We look now at any longer  bispecial $w_t$ resolved by rule $(c)$:  as there are only rules $(c)$ and $(x)$ (strictly)  between  $\lambda_2$ and $\mu_y$, the  $yw_{q'}y$ defined above is
 $u$ as long as  $\lambda_2\leq t\leq \mu_y$, and we know that $u$ does not occur in $Z_p^2$.  Thus all these bispecials satisfy the order condition in $\Lambda_{Z_p}$; as for still longer bispecial words of $\Lambda$, they are too long to occur  in $Z_p^2$, as $w_{\mu_y+1}=w_{\mu_y}L_{\mu_y}$ has a length greater than $2|L_{\mu_y}|$.

There remain to consider the $A_p$, $B_p$ or $C_p$ for initial values of $p$. It is immediate that those of the form $ca^k$ or $ba^k$ do cluster, while the $C_p$, $\lambda_1+1\leq p\leq \lambda_2-1$, are dealt with as in Proposition \ref{sq}, by changing the language and checking that $p\leq \mu'_y$ in the new language, and all these cluster. \qed\\

The following statements give an equivalent criterion for $A_p$, $B_p$ and $C_p$ to cluster, which gives more information and does not particularize any order of apparition of the rules. 

\begin{corollary}\label{clist} 
With or without the assumption $(Rabc)$, $Z_p$ clusters if and only if there exists a letter $z'$ in $\{a,b,c\}$ such that in the word $D_pz$
\begin{itemize}
\item if $z'=a$, neither the letters $z'$, $b$,  $c$ nor the letters $z'$, $c$,  $b$, 
\item if $z'=b$,   neither the letters $z'$, $a$,  $c$ nor the letters $z'$, $c$,  $a$, 
\item if $z'=c$, neither the letters $z'$, $b$,  $a$ nor the letters $z'$, $a$,  $b$, 
\end{itemize}
occur at any increasing sequence of indices.\\
Moreover, in the cases where $Z_p$ clusters, it does cluster perfectly for any order for which a possible  $z'$ is in the middle of $\{a,b,c\}$, and does not cluster for any other order and permutation. 
\end{corollary}
{\bf Proof}\\
 Suppose first $(Rabc)$ holds. Then all the assertions but the last two are deduced directly from Proposition \ref{list}. As for the last two, they are consequences of  Proposition \ref{list} when $z'=c$, with orders dictated by Lemmas \ref{caro} and \ref{rec}, and are proved in the same way for $z'=a$ or $z'=b$. The other cases for the order of apparition of the rules are deduced by a permutation of the letters, under which the conclusion is invariant. \qed\\

 Thus, as clustering is invariant by conjugacy, we know all the clustering AR words satisfying the assertions of Proposition \ref{sq}. We know also all the clustering AR words which are conjugate to a power of a standard AR word, or equivalently can be de-substituted to a power of one letter by way of the six AR morphisms, they are the powers of the standard  clustering words of  Proposition \ref{list} or Corollary \ref{clist} and  all their conjugates.
 
 But this does not tell which ones are in a fixed language. For Tribonacci, and all AR languages with the same first five rules, the  $A_p$, $B_p$ or $C_p$ which cluster (perfectly) for the order $a<c<b$ (or its symmetric) are $A_0$ to $A_4$, $B_0$ to $B_3$, $C_0$ to  $C_4=cabaabacaba$, of length $11$, the longest standard word which clusters, while $B_4=bacabaabacaba$, of length $13$, is the shortest standard one which does not cluster. Also, $A_0$ to $A_2$, $B_0$ to, $B_3$, $C_0$ to $C_3$ cluster (perfectly) for the order $a<b<c$ (or its symmetric),  $A_0$ to $A_2$, $B_0$ to $B_2$, $C_0$ and $C_1$ cluster (perfectly) for the order $b<a<c$ (or its symmetric). Other clustering words conjugate to standard or  powers of standard  are for example the non-primitive $A_3^3=(abacaba)^3$, of length $21$, or  $(abacaba)^2cabaabacaba$  of  Theorem \ref{arc} or Proposition \ref{long} below;  for these last two words $v$,  $vv$ is not in the Tribonacci language. \\

In  a given AR language $\Lambda$, we find now a clustering word, conjugate to a standard AR word,  which is longer than all the ones in Proposition \ref{list}.

\begin{proposition}\label{long} Let $\Lambda$ be an AR language satisfying $(Rabc)$. With the notations of Theorem \ref{arc}, let also $\mu$ be the stage of the first rule in the string of rules $(z)$, $z=x$ or $z=c$,  just  before stage $\mu_y$. The word
$v=S_{\mu_y}^{\mu_y-\mu+1}M_{\mu_y}$ is a primitive perfectly clustering (for the order $a<c<b$ or its symmetric) word of $\Lambda$ conjugate to a standard AR word. .\end{proposition}
{\bf Proof}\\
We have $S_{\mu_y}=Z_{\mu_y}$ and, by Lemma \ref{lms}, $L_{\mu_y}=Y_{\mu_y}$. We define another AR language $\hat\Lambda$ by the same rules as $\Lambda$ up to (and including) stage $\mu_y-1$, then $\mu_y-\mu+1$ rules $(z)$, and any acceptable sequence of rules beyond. Then $\hat\Lambda$ has the same $x$ and $y$ as $\Lambda$, and $v$ is conjugate to $\hat X_{2\mu_y+\mu-1}$ ( (of $\hat\Lambda$) if $z=c$,  to $\hat C_{2\mu_y+\mu-1}$  if $z=x$, and $2\mu_y+\mu-1<\hat\mu_y$ (of $\hat\Lambda$), thus by Proposition \ref{list} $v$ is a primitive  perfectly clustering word. \\

It remains to prove that $v$ is in $\Lambda$. We
know that $w_{\mu_y+1}L_{\mu_y+1}=w_{\mu_y}L_{\mu_y}M_{\mu_y}L_{\mu_y}$ is in $\Lambda$. We have $\mu\geq \mu_x$; suppose first that either $\mu>\mu_x$, or $\mu=\mu_x$ and $\lambda_y<\lambda_x$. Then by Lemma \ref{lms}  $Z_{\mu}=M_{\mu}$. Thus $S_{\mu_y}=Z_{\mu}$ while $L_{\mu_y}=L_{\mu}Z_{\mu}^{\mu_y-\mu}$,
 $M_{\mu_y}=M_{\mu}Z_{\mu}^{\mu_y-\mu}$.
Then $v=M_{\mu}^{\mu_y-\mu+1}S_{\mu}M_{\mu}^{\mu_y-\mu}$, while
$L_{\mu_y}M_{\mu_y}=L_{\mu}M_{\mu}^{\mu_y-\mu}S_{\mu}M_{\mu}^{\mu_y-\mu}$ is in $\Lambda$, and all we have to prove is that $M_{\mu}$ is a suffix of $L_{\mu}$, which is true by Lemma \ref{rel} as $\mu \geq \lambda_2
+1$.\\

Suppose now that  $\mu=\mu_x$ and $\lambda_x<\lambda_y$. 
Then $z=x$, $Z_{\mu}=X_{\mu}=L_{\mu}$, $v=L_{\mu}^{\mu_y-\mu+1}S_{\mu}L_{\mu}^{\mu_y-\mu}$, while
$w_{\mu_y}L_{\mu_y}M_{\mu_y}=w_{\mu_y}M_{\mu}L_{\mu}^{\mu_y-\mu}S_{\mu}L_{\mu}^{\mu_y-\mu}$ is in $\Lambda$, thus what we have to prove is that $L_{\mu}=X_{\mu}$ is a suffix of $w_{\mu_y}M_{\mu}=w_{\mu_y}Y _{\mu}$, which will be true if $X_{\mu}$ is a suffix of $w_{\mu}Y_{\mu}$, as  $w_{\mu_y}$ ends with  $w_{\mu}$. 

Going backward through rules $(c)$, what we have to prove is that $X_{\lambda_2}$ is a suffix of $w_{\mu}Y_{\lambda_2}$. Then the rule at stage $\lambda_2-1$ is $(y)$, thus we have to prove that $X_{\lambda_2-1}$ is a suffix of $w_{\mu}$, and this is true as  $w_{\lambda_2-1}$  is a suffix of  $w_{\mu}$, and  by Lemma \ref{rel} $X_{\lambda_2-1}$ is a suffix of   $w_{\lambda_2-1}$, except possibly if $x=b$ and $\lambda_2-1=\lambda_1$, which cannot happen as then the  rule at stage $\lambda_1=\lambda_2-1$ should be both $(b)$ and  $(y)=(a)$. 
\qed\\

\begin{corollary}\label{gen} For every $n$, there  exist arbitrarily long primitive (perfectly) clustering Arnoux-Rauzy words with at least $n$ occurrences of each letter.\\
Every Arnoux-Rauzy language contains a primitive (perfectly) clustering word of length at least $22$.\end{corollary}
{\bf Proof}\\ 
Take an AR language where $D$ begins with $abc^n$. Then we get  $B_p=ba(caba)^n$ for some $p\leq \mu_x$.

For a general AR language, the smallest possible value of the length of the word in Proposition \ref{long} is $|S_{\mu_a}S_{\mu_a}M_{\mu _a}|$  where $D$ begins with $abcba$, which gives $22$. \qed\\

\begin{corollary}\label{sha}  When the directive word $D$ begins with $ab^{n_1}c^{n_2}a^{n_3}b$
 for any integers $n_i\geq 1$, $i=1,2,3$,  the word of Proposition \ref{long} has maximal length among clustering words of $\Lambda$, and the bound in Theorem \ref{arc} is optimal. Assuming $(Rabc)$, in all other cases, there is  a gap between the length of the word in Proposition \ref{long} and the bound in Theorem \ref{arc}.\end{corollary}
{\bf Proof}\\
We look at he  proof of Proposition \ref{long}. In the case where either $\mu>\mu_x$, or $\mu=\mu_x$ and $\lambda_y<\lambda_x$, there is always a gap of at least $2$ between
the best bound in Theorem \ref{arc} and the length of the word in Proposition \ref{long}.  In the case where  $\mu=\mu_x$ and $\lambda_x<\lambda_y$, this gap is reduced to $1$ whenever  $X_{\lambda_2}=w_{\lambda_2-1}Y_{\lambda_2}$, which is equivalent to $D$ being as in the assertion above. \qed\\

\begin{example}\label{e1}\rm Take the Tribonacci language or any AR language where the directive word $D$ begins with $abcab$. The word in Proposition \ref{long} is the example of length $25$ in 
 Theorem \ref{arc}, which thus is a  conjugate to a standard AR word. It is a palindrome, and we check that it is the only  clustering word of maximal length in $\Lambda$.\end{example}
 
\begin{example}\label{e2}\rm Take any 
AR language $\Lambda$ where $D$ begins with  $abacba$, the bound in Theorem \ref{arc} is $45$, the  word in Proposition \ref{long} has length $43$. We have seen  in the proof of Theorem \ref{arc} that the only word of length $44$ which might cluster  is $v=abaaba(cabaababaabacabaaba)^2$ and that it does not cluster. Indeed $vv$  does contain all the words $zw_{\lambda_1}z'$ and $zz'$ of $\Lambda$, $z\neq z'$,  and we check that for $w_{\lambda_2}=abaaba$, the four  extensions $aw_{\lambda_2}c$, $bw_{\lambda_2}c$, $cw_{\lambda_2}a$, $cw_{\lambda_2}b$  appear in $vv$ as do indeed $cw_{\lambda_2}c$ and $aw_{\lambda_2}a$, thus $v$ does not cluster and  $v$ is not  conjugate to a standard AR word. Thus  the  word in Proposition \ref{long}, which is $v$ deprived of its first letter, is a clustering word of maximal length  in $\Lambda$; it is not the only one, as its reverse is also clustering.\end{example}

 \begin{example}\label{e3}\rm Take any 
AR language $\Lambda$ where $D$ begins with  $abcba$,  the bound in Theorem \ref{arc} is $24$, the  word in Proposition \ref{long} has length $22$. We have seen  in the proof of Theorem \ref{arc} that the only word of length $23$ which might cluster  is $v'=a(bacaba)^2cababacaba$ and that it does not cluster. Indeed $v'$ is $L_{\mu_a}M_{\mu_a}$ thus is conjugate to $L_{\mu_a+1}$, hence by  Proposition \ref{sq} $v'$ is  conjugate to a standard AR word, but does not cluster. Thus  the  word in Proposition \ref{long}, which is $v$ deprived of its first letter, is a clustering word of maximal length  in $\Lambda$;  it is not the only one, as its reverse is also clustering.\end{example}

\begin{conj}\label{clong}  In a given AR language $\Lambda$ satisfying $(Rabc)$,  the word in Proposition \ref{long} is the longest clustering word,  or, if this fails,  the longest clustering word  conjugate to a standard AR word. \end{conj}

\subsection{Clustering AR words non conjugate to standard}\label{sns}
 
  We turn now to  words which are not  conjugate to  standard AR words.

  \begin{example}\label{e5}\rm 
 For all $n$,  $ba(ca)^nb$ is  an AR word not  conjugate to a standard AR word, and it is perfectly clustering.  \end{example}
  
The following propositions characterize, in two steps,  all the words having this property, by identifying the particular way they are generated in  the general construction of  \cite{sp}.

\begin{proposition}\label{ptb1} Let $w$  on the alphabet $\mathcal A=\{a,b,c\}$ be a perfectly clustering AR word which is not in the range of any of the six AR morphisms $\tau_x$ and $\sigma_x$ for $x\in \{a,b,c\}.$ Then, up to a permutation of the letters $a,b$ and $c,$ there exists a word $v$ on the alphabet $\{a,c\}$ containing both $a$ and $c$ such that the conjugate $w'=b^{-1}wb$  is obtained from $\tau_v(b)$ by inserting a single $b$ between each pair of consecutive occurrences of $a$ or between each pair of consecutive occurrences of  $c$  in $\tau_v(b)$ (where at most one of $aa$ and $cc$ can occur) plus a $b$ at the very beginning.  Furthermore, $w$  is a palindrome beginning and ending in $b$ containing both $a$ and $c$ but no $a^2$ nor $c^2$, is primitive, and  any order for which  $w$ is perfectly clustering has $b$ as the middle letter. 

Conversely, any word $w$ built as above is a perfectly clustering AR word for the order $a<b<c$, not in the range of any of the six AR morphisms.
\end{proposition}

{\bf Proof}\\
Let $w$ be as in the first sentence.  Then $w$ is not a power of a single letter and $|w|\geq 5.$ Let $a$ denote the separating letter of $w.$ Since $w$ is not in the range of $\tau_a$ nor $\sigma_a,$ $w$ begins and ends in some letter different from $a.$ Let $b$ denote the first letter of $w.$ Then $w$ also ends in $b$, because otherwise $cb$ is in $\Lambda_w$ but not $bc$, while $\Lambda_w$ is closed under reversal by Proposition \ref{rev}. Also, since $a$ is the separating letter of $w$, $bb$ does not occur in $w$, although it occurs in $\Lambda_w$. Thus $w$ is a palindrome,  as $\bar w$ must be  conjugate to $w$, but the only conjugate of $w$ which does not contain $bb$ is $w$ itself. Also, $w$ is  primitive: if $w=v^n$ for some $n\geq 2,$ then as $v$ must begin and end in $b,$ $bb$ is a factor of $w,$ a contradiction.\\

{\bf Claim 1 :} Each letter of $\mathcal A$ must occur in $w.$

By assumption each of $a$ and $b$ occurs in $w.$ If $c$ does not occur in $w,$ then $w$ is a perfectly clustering binary palindrome of the form $w=bub$ where $u$ begins and ends in the letter $a.$ Furthermore $w$ cannot contain a factor of the form $ba^nb$ with $n> 1$, for otherwise $\Lambda_w$ contains both $aa$ and $bb,$ which contradicts Theorem \ref{car}. Moreover, as $bb$ does not occur in $w,$ it follows that any two consecutive occurrences of $b$ in $w$ must be separated by a single $a$. Thus $w=b a(ba)^nb$ for some $n\geq 0$ and hence $w=\tau_b(a^{n+1}b)$ contradicting that $w$ is not in the range of $\tau_b.$\\

As $w$ is perfectly clustering and $bb\in \Lambda_w,$ by Theorem \ref{car} $b$ must be the middle letter under any (perfectly) clustering order on $A$, and furthermore $aa$ does not occur in $w.$  Let us consider the conjugate $w'=b^{-1}wb.$ Note that $w'$ begins in $a$ and ends in $abb.$ \\

{\bf Claim 2 : }$w'=\psi (\tau_v(b))$ for some word $v$ on the alphabet $\{a,c\}$ beginning in $a$ and containing $c$ and where the mapping $\psi$, defined in Section 3 of \cite{sp}, amounts to inserting a single $b$ in the middle of each occurrence of $aa$  and $ab$ in $\tau_v(b).$ 
  $\tau_u(w)$ is a perfectly clustering AR word for the order $a<b<c.$ Furthermore, every conjugate of $\tau_u(w)$ different from $\tau_u(w)$ is not an AR word and hence in particular,  $\tau_u(w)$ is not conjugate to a standard AR word.
We note that $w'$ is a perfectly clustering word (for the order $a<b<c)$ and no conjugate of $w'$ is in the range of $\tau_a$ nor $\tau_c.$ In fact,
every conjugate of $w'$ (other than $w)$ contains $bb$ as a factor, hence is not in the range of $\tau_a$ nor $\tau_c,$ and by assumption the same is true for $w.$  By application of Lemmas~3.7 and 3.8 of \cite{sp},  $w'=\psi(ub)$ where $u$ is a word on the alphabet $\{a,c\}$ and $ub$ is also perfectly clustering relative to the order $a<b<c.$ Note that if $b$ occurred in$u$, then $w$ would admit an occurrence of $bb$ contrary to our assumption that $a$ is the separating letter of $w.$ Thus $ub$ contains both $a$ and $c$ and $\Lambda_{ub}$ contains each of $ab,ba,ac,ca.$ It follows from Lemma 3.1 of \cite{sp} that $cc$ is not in  $\Lambda_{ub}$ and hence each occurrence of $c$ in $ub$ must be directly preceded and followed by the letter $a.$ In other words, $a$ is a separating letter of $ub$ and $w'$ is obtained from $ub$ by inserting a single $b$ in the middle of each $aa$ occurring in $u,$ plus an additional $b$ at the end. Also $u$ is a palindrome. 

By an iterated application of Lemma 3.4 in \cite{sp}, we can write $ub=\tau_a^r(u'b)$ for some $r\geq 1$ and $u'$ on the alphabet $\{a,c\}$ beginning in $c$ and furthermore $u'b$ is perfectly clustering for $a<b<c.$ Thus $u'$ is also a palindrome and hence $bc$ and $cb$ are both in $\Lambda_{u'b}.$ If the letter $a$ does not occur in $u',$ then we can write $u'b=c^sb=\tau_{c^s}(b)$ for some $s\geq 1$ and hence \[w'=\psi(ub)=\psi (\tau_{a^r}(u'b))=\psi (\tau_{a^r}\circ \tau_{c^s}(b)=\psi (\tau_{a^rc^s}(b))\] as required. On the other hand, if $a$ occurs in $u',$ then it follows from Lemma 3.1 in \cite{sp} that each $a$ in $u'$ must be preceded and followed by the letter $c$ and thus $c$ is a separating letter of $u'b.$ Thus similarly we can write $u'b=\tau_c(u''b)$ where $u''b$ is perfectly clustering for the order $a<b<c$ and $u''$ is a palindrome beginning and ending in $a$ or $c.$ Continuing in this way, we eventually obtain that $ub=\tau_v(b)$ for some word $v$ on the alphabet $\{a,c\}$ containing each of $a$ and $c.$ \\

In the other direction, let $v$ be a word on the alphabet $\{a,c\}$ containing each of $a$ and $c.$ Without loss of generality, we may assume that $v$ begins in the letter $a.$ It follows from Lemmas~3.3 and 3.4  of \cite{sp} that $\tau_v(b)$ is perfectly clustering for the order $a<b<c.$ Also clearly $\tau_v(b)$ is an AR word. By definition, $w'=\psi (\tau_v(b))$
where $\psi$ is the mapping defined in Section 3 of \cite{sp}. It follows from Lemma~3.7 of \cite{sp} that $w',$ and hence $w,$ is perfectly clustering for the order $a<b<c.$ It remains to show that $w$ is an AR word which is not in the range of any of the six AR morphisms. We first note that by definition $w$ begins and ends in $b$ and hence can only be in the range of $\tau_b$ or $\sigma_b.$ But as $w$ contains $ac$ neither is possible. Finally it is easily verified that $b^{-1}wb^{-1}=\tau_{av'}(a)$ where $v'$ is obtained from $a^{-1}v$ by exchanging the letters $a$ and $b$ and keeping $c$ fixed. It follows that
$b^{-1}wb^{-1}$ is a bispecial AR word from which it follows that $w$ is an AR word as required.  \qed\\

We note that the shortest $w$ verifying the assumptions of Proposition \ref{ptb1} is (up to a permutation of the letters) $bacab$, built from $\tau_{ac}b$.

\begin{proposition}\label{ptb2} Assume $v$ is a perfectly clustering AR word which is not conjugate to a standard AR word, nor to any power of a standard AR word. Then up to a permutation of the letters, $v$ is conjugate to a word of the form $\tau_u(w)$ where $u$ (possibly empty) is on the alphabet $\{a,c\}$ and where $w$  
is as in Proposition~\ref{ptb1}.

  Conversely, let $v$ be as in the previous sentence; then it is an AR word perfectly clustering for the order $a<b<c.$ Furthermore, every conjugate of $v$ different from $v$ is not an AR word and hence $v$ is not conjugate to any power of a  standard AR word. \end{proposition}

{\bf Proof}\\ 
By Proposition \ref{sq}, an AR word is conjugate to a power of a standard if and only if it can be de-substituted to a power of one letter by way of the six AR morphisms. This means that otherwise we can write $v=f(w)$ where $f$ (possibly the identity)  is some concatenation of  AR morphisms, and where $v$ is not in the range of any of the six AR morphisms. Assume $w$ is perfectly clustering for the order $a<b<c,$ then $f$ cannot involve $\tau_b$ nor $\sigma_b$ and hence is a concatenation of  $\tau_a$, $\tau_c$, $\sigma_a$,  $\sigma_c\}$; replacing $v$ by a conjugate, we can get that $f$ is a concatenation of  $\tau_a$, $\tau_c$.  By  Lemmas~3.3 and 3.4 in \cite{sp}, we have that $v$ is also perfectly clustering for $a<b<c.$ \\

In the other direction, we begin by considering the case when $u$ is empty. It follows from Proposition~\ref{ptb1} that $v$ is perfectly clustering AR word for the order $a<b<c.$ 
Furthermore $v$ contains each of $a,b$ and $c,$ begins and ends in $b$ and has either $a$ or $c$ as a separating letter.Now let $v'$ be a conjugate of $v$ with $v'\neq v.$ Then $v'$ contains $bb$ and either $ac$ or $ca$ or both and hence $v'$ is not an AR word.

Now assume that $u$ is a non-empty word on the alphabet $\{a,c\}.$ The properties of $v$ being an AR word containing each letter, and no conjugate $v'\neq v$ being an AR word, are clearly stable under application of $\tau_u.$ By Lemmas 3.3 and 3.4 in \cite{sp}, the clustering property is also stable by application of $\tau_u.$   
 \qed\\

However, there are also  infinitely many primitive  AR words which cluster but  not perfectly, and hence are not conjugate to a standard AR word.

    \begin{example}\label{e4}\rm  The word $w=abaca$, which belongs to the Tribonacci language, does cluster for the order $a<c<b$ (and no other one), for the permutation $\pi a=c$, $\pi b=a$, $\pi c=b$ (thus not perfectly), and in the language $\Lambda_w$ the bispecial $a$ is resolved by $caa$, $aab$, $bac$, thus the bispecials of $\Lambda_w$ are not in any  AR language, though they all satisfy the order condition. Thus $abaca$ is not conjugate to a standard AR word, and it is  clustering but not perfectly clustering. The same properties are shared by  $aba^nca^n$ for all $n$. \end{example}

 \begin{question}\label{q} What are the  primitive  AR words which cluster but not perfectly? \end{question}

As for general clustering words on three letters, they are characterized in \cite{sp} for the symmetric permutation, and in \cite{fz4} for all permutations, and are not always AR words.

\begin{proposition} For any order and any  permutation $\pi$  different from the identity, there are infinitely many $\pi$-clustering words on three letters which are not AR words. \end{proposition}
{\bf Proof}\\
We fix an order and a permutation. By \cite{fz4}, any word $w$ such that $ww$ is in the language $\Lambda$ corresponding to  an interval exchange transformation built with this order and this permutation, and satisfying the {\em i.d.o.c. condition}, is clustering for this order and this permutation. Such a $\Lambda$ is uniformly recurrent, and contains infinitely many squares, by \cite{fz1} if $\pi$ is the symmetric permutation, \cite{fie} in  general, thus infinitely many such clustering words $w$. If a word $w$ in $\Lambda$ contains all the extensions $xvy$ of all bispecial words $v$ in $\Lambda$ longer than some constant,  $w$ cannot be  an AR word by \cite{fz3}. And this will be true for any long enough clustering word $w$ in $\Lambda$. \qed

\section{Generalizations}\label{sgen}

\subsection{Sturmians}\label{sstu}
On two letters $a$ and $b$, as is shown in \cite{ar}, the Sturmian languages  of \cite{hm} can be generated by words $A_k$ and $B_k$,  which are called {\em standard} Sturmian words, using AR-type rules on two letters. Each Sturmian language contains infinitely many clustering words, and all these are known since \cite{mrs} and \cite{jz}. For sake of completeness, we reprove this result by the methods of the present paper.

\begin{proposition}\label{tstu} The primitive clustering Sturmian words, as well as the primitive clustering words on $\{a,b\}$, are  all the standard Sturmian words and all their conjugates. \end{proposition}
{\bf Proof}\\
Note first that the only clustering words with the identity as permutation are the $a^m$ and $b^m$, thus we can restrict ourself to perfectly clustering words, for the order $a<b$. The Sturmian languages are identified in \cite{hm} with $2$-interval exchange languages, thus we deduce from \cite{fz4} that a primitive Sturmian word $v$, or a primitive word $v$ on $\{a,b\}$, is clustering iff $vv$ is a factor of  a Sturmian language. This, by the same proof as  Proposition \ref{sq}, is equivalent to $v$  conjugate to some $A_p$ or $B_p$ in some Sturmian language. \qed\\

Note that, for Sturmian languages or more generally for interval exchange languages, the necessary and sufficient  condition for $v$  to cluster in Theorem \ref{car} is, by \cite{fhuz}, equivalent to the one given in \cite{fz4}, namely that $vv$ be  a factor of such a language. 
However, to determine if a given word clusters, our Theorem  \ref{car} is more explicit. Take for example the two Sturmian words $v=abaa$ and $v'=baab$, $v$ and $v'$ are factors of the {\it Fibonacci} language, while $vv$ and $v'v'$ are not in this language. It is easy to check by hand that $v$ clusters (for $a<b$ and the symmetric permutation) and $v'$ does not cluster for any order and permutation; thus we know that $vv$ must be in some Sturmian language and $v'v'$ cannot be in any Sturmian language, but it is easier, and quicker in general than computing the Burrows-Wheeler transform, to  check directly that the bispecials in $vv$ satisfy the order condition and those in $v'v'$ do not; in this last case, it is immediate that no order condition can be  satisfied by the empty bispecial, as $aa$, $bb$, $ab$, $ba$ occur in $v'v'$.

\subsection{Episturmians on three letters}\label{sepi}
In the literature, for which we refer the reader to the two surveys \cite{be} and \cite{gjp}, we found only the definition of episturmian infinite words, one-sided in general though two-sided words are briefly considered in \cite{gjp}. To make the  present paper coherent, we define here the corresponding languages, our definition having been chosen to correspond to what is used in practice.

\begin{definition} A language is {\em episturmian} if it is uniformly recurrent, closed under reversal, and admits at most one right special factor of each length.\end{definition}

An episturmian language on three letters can be generated by AR rules or by AR morphisms, as is proved in the founding paper \cite{djp}. Indeed, these episturmian languages can be defined by a modification of Definition \ref{dar} above, where the assumption ``each one of the three rules is used infinitely many times" is replaced by ``each one of the three rules is used at least once".

The description of the bispecial words is deduced from the one given after Definition \ref{dar} by the following modifications: the possible bispecials are among the $w_k$, and $w_k$ has at most three suffix return words which are among $A_k$, $B_k$ and $C_k$. More precisely, $A_k$ is a return word of $w_k$, or equivalently the label of an elementary circuit in the Rauzy graphs, if and only if $w_kA_k$, or equivalently $w_ka$, is in $\Lambda$, and similarly for $B_k$ and $C_k$. 

\begin{lemma}\label{ebs} The word $w_pA_p$ is in an episturmian language $\Lambda$ on three letters if and only if the directive word of $\Lambda$ is such that there exist rules $(a)$ at or after stage $p$, and similarly for $B_p$ and rules $(b)$, $C_p$ and rules $(c)$. 
\end{lemma}
{\bf Proof}\\ Suppose for example there is a rule $(a)$ at or after stage $p$. Then, for some $q\geq p$, $w_{q+1}=w_qA_q$, thus $w_qA_q$ is in $\Lambda$, and so is $w_pA_p$ as $w_p$ is a suffix of $w_q$ and  $A_p$ is a prefix of $A_q$.

Suppose  there is no rule $(a)$ at or after stage $p$. Then $w_pA_p$ cannot be in $B_r$ or $C_r$ as $B_r$ and $C_r$ do not have
 $A_p$ in their decomposition by AR rules. As there are infinitely many rules $(b)$ or $(c)$, the length of $B_r$ or $C_r$ tends to infinity with $r$, thus this contradicts uniform recurrence.  \qed\\                               

The assertion $(Rabc)$ is defined as for AR languages. 

\begin{theorem}\label{thepi} An episturmian language $\Lambda$ on three letters, satisfying $(Rabc)$, contains infinitely many clustering words if and only if its directive word is $D'D"$, where $D'$ is a finite word on the alphabet $\{a,b\}$ and $D"$ is a one-sided infinite word on the alphabet $\{a,c\}$ or $\{b,c\}$.\\
When this is not the case, with the notations of Theorem \ref{arc}, no word of length at least\\
 $|w_{\lambda_2}|+\max\{|Z_{\mu_y+1}|, Z \in  \{A,B,C\}, w_{\mu_y+1}Z_{\mu_y+1}\in \Lambda\}+1$\\ can cluster for any permutation  and any order.\end{theorem}
{\bf Proof}\\
Suppose $D$ is $D'D"$ as in the hypothesis. Let $v$ be any word such that $vv$ is in $\Lambda$. The bispecials in the language $\Lambda_v$ are some of the $w_p$ of $\Lambda$, and are resolved as in the AR rules, possibly with some extensions missing. By Lemma \ref{rec}, those which are resolved as in  rule $(a)$ or $(b)$ satisfy the order condition for $a<c<b$. For those $w_p$ which are resolved as in rule $(c)$,  there is no rule $(z)$ on or after stage $p$ for $z=a$ or $z=b$, thus by Lemma \ref{ebs} $w_pz$ is not in $\Lambda$, and $w_p$ is resolved by a subset of $\{aw_pc, cw_pa, cw_pc\}$ or $\{bw_pc, cw_pb, cw_pc\}$, thus satisfies the order condition. Thus by Lemma \ref{caro} $v$ clusters perfectly. To find such $v$, we follow the reasoning of Proposition  \ref{sq}: any $Z_p$ in  $\{A_p,B_p,C_p\}$ will have its square in $\Lambda$ provided $p$ is large enough and $w_pZ_p$ is in $\Lambda$. If each letter in $D"$ occurs infinitely many times, we get arbitrarily long  primitive  clustering words; otherwise all the $Z_p^n$ will be clustering for one value of $Z$. 

Suppose $D$ is not $D'D"$ as in the hypothesis. Then after the first rule $(c)$ there is at least one rule $(a)$ and one rule
$(a)$. Thus all the quantities in Theorem \ref{arc} can be defined, and we can follow the reasoning of this theorem, with the following modifications: the assertion $(OA_q)$ is now that  $v$ occurs once in $w_qA_q$,  does not occur in $w_qC_q$ or $w_qB_q$, and $w_qA_q$ is in $\Lambda$, and similarly for $(OB_q)$ and $(OC_q)$; the assertion $(OT_q)$ is now that $v$ occurs at least once in each $w_qZ_q$ which is in $\Lambda$, $Z$ in  $\{A,B,C\}$, and the maximal return time of $v$ is the maximal length of these $Z_q$. Then, we get the maximal return times of the two special words, and conclude immediately as, in contrast with Theorem \ref{arc}, we keep the quantity $|w_{\lambda_2}|$ in the bound.\qed\\

 \begin{example}\label{e6}\rm Let $D=abc(ab)^{\omega}$. This  gives an episturmian language which contains only finitely many clustering words, but is not an AR language. Its complexity function is $p(n)=2n+1$ for $1\leq n\leq 4$, $p(n)=n+5$ for $n\geq 5$. \end{example}
 
Note that these properties are shared by all episturmian languages where $D=abcD"$
where $D"$ is a one-sided infinite word on $\{a,b\}$ in which both $a$ and $b$ occur infinitely many times. By Theorem \ref{thepi}, any episturmian language whose complexity is at least $n+1$ for all $n$ but is strictly smaller than the $p(n)$ of Example \ref{e6} produces infinitely many clustering words. One can wonder whether this is still  true for  any  language, or at least for  any uniformly recurrent language,  of complexity  at least $n+1$ for all $n$ but  strictly smaller than this $p(n)$. For sake of completeness, we give a (non episturmian) counter-example.

\begin{example}\label{e7}\rm We build a language $\Lambda$ on $\{a,b,c\}$ in the following way: the empty bispecial word is resolved by $\{ab,ac,ba,ca\}$; the bispecial $a$ is resolved by $\{bab,bac,cab,cac\}$;  the bispecial $aba$ is resolved by $\{babab,babac,cabab\}$;  the bispecial $aca$ is resolved by $\{bacac,cacab\}$; every further bispecial $w$ is resolved either by $\{bwb,bwc,cwb\}$ or by $\{bwc, cwb, cwc\}$, each possibility being used infinitely many times.  Its complexity function is $p(1)=3$, $p(2)=4$, $p(3)=6$, $p(n)=n+4$ for $n\geq 4$. Its Rauzy graphs of length $4$ and more have the same shape as the Rauzy graphs of Sturmian languages, thus the alternating of resolution rules ensures that $\Lambda$ is uniformly recurrent. But the bispecial word $a$ does not satisfy the requirement of the order condition, for any order and permutation,  and its four extensions $bab,bac,cab,cac$ occur in every long enough factor of $\Lambda$. Thus by Theorem \ref{car}  $\Lambda$ contains only finitely many clustering words.\end{example}

\subsection{Larger alphabets}\label{slga}
AR languages can be generalized to any alphabet
$\mathcal A=\{a_1,a_2,\cdots, a_r\}$ (note that here the order will not  necessarily be $a_1<a_2< \cdots <a_r$). 

\begin{definition} An AR language is generated by words $A_k^{(i)}$, $1\leq i\leq r$, starting from $A_0^{(i)}=a_i$, $1\leq i\leq r$, and by rule $(a_i)$
 at stage $k$, $A^{(i)}_{k+1}=A^{(i)}_{k}$,
$A^{(j)}_{k+1}=A^{(j)}_kA^{(i)}_k$ for all $i\neq j$.

Each rule is used infinitely many times. The directive word is defined in the usual way. The $A^{(i)}_k$ are again the labels of the elementary circuits in the Rauzy graphs.
The {\em $r$-Bonacci language}, $r\geq 3$, is defined by $D=(a_1...a_r)^{\omega}$. \end{definition}

There the methods of Section \ref{sar} apply mutatis mutandis, but the number of cases to be considered grows very quickly, and a lot of space would  be required to generalize all the above study. Thus we shall just generalize Theorem \ref{arc}, with some  loss of optimality. \\

\begin{proposition}
We denote by $(a)$, $(b)$, $(c)$ the first three rules by order of appearance, and define all quantities in Theorem \ref{arc} above in the same way. Let $v$ be an AR word on an $r$-letter alphabet of length at least $|w_{\lambda_2}|+\max_{1\leq i\leq r}(|A^{(i)}_{\mu_y+1}|)+1$. Then $v$ cannot cluster for any permutation  and any order on the alphabet.\\
For $4$-Bonacci, this bound is not optimal; the better bound $|w_{\lambda_y}|+\max_{1\leq i\leq r}(|A^{(i)}_{\mu_y+1}|)+1$ holds but is not optimal either.
\end{proposition}
{\bf Proof}\\
 Then Lemma \ref{caro} is still valid, by restricting the orders on $\mathcal A$ to the set $\{a,b,c\}$, and again we need to know the maximal return times of 
$u'=xw_{\lambda_x}x$ and $u=yw_{\lambda_y}y$. These are computed exactly as in Theorem \ref{arc}, mutatis mutandis:  the assertions are now $(O^{(j)}_q)$, that $v$ occurs once in $w_qA^{(j)}_q$ and does not occur in any $w_qA^{(i)}_q$, $i\neq j$ ($a,b,c$ being denoted also by $a_1,a_2, a_3$), and $OT_q$, that $v$ occurs at least once in each  $w_qA^{(i)}_q$ and the maximal return time of $v$ is $\max_{1\leq i\leq r}(|A^{(i)}_q|)$. These evolve like in the proof of Theorem \ref{arc}, and thus the maximal return time of $u$ is $\max_{1\leq i\leq r}(|A^{(i)}_{\mu_y}|)$, the maximal return time of $u'$ is $\max_{1\leq i\leq r}(|A^{(i)}_{\mu_x})|$. We conclude immediately as we keep the quantity $|w_{\lambda_2}|$ in the bound.

 For $4$-Bonacci, if we denote the letters by $a$,  $b$, $c$, $d$, we have $\lambda_2=2$, $w_3=abacaba$, $w_2=aba$, $w_1=a$, $w_0$ is the empty word. The bound in the conclusion is $|C_6|+4=60$, but in this simple case we can mimic the end of the proof of Theorem \ref{arc} and  replace the $|w_{\lambda_2}|$ in the bound by $|w_{\lambda_y}|$, thus getting 
$|C_6|+2=58$. As $C_6=C_5B_5$, by the usual reasoning the only word of length  $|C_6|+1$ which does not contain $bw_{\lambda_2}c$ nor $cw_{\lambda_2}b$ is, up to cyclic conjugacy,  $u=aB_5C_5$. In the language $\Lambda_u$, we check that $w_0$ is resolved by $\{aa, ab,ac,ad, ba,ca,da\}$,  $w_2$ is resolved by $\{cw_2c, aw_2c, cw_2a,  dw_2c, cw_2d\}$,  $w_3$ is resolved by $\{aw_3d, dw_3a,  dw_3c, cw_3d\}$, thus, if $u$ clusters, by Lemma \ref{ord} each one of  $a$, $c$, $d$ must be at an end of the order between them, thus no word of length $57$ can cluster for any order and permutation. \qed\\

In the general case, we do not try to replace the $|w_{\lambda_2}|$ in the bound by $|w_{\lambda_y}|$ as the proof would be complicated by the presence of rules $(a_i)$, $i\geq 4$, between stages $\lambda_2$ and $\mu_y$, and the improved bound is not optimal even for $4$-Bonacci.\\

Similarly, the main result of Section \ref{slga} can be generalized to episturmian languages on larger alphabets: an episturmian language $\Lambda$ on $r$ letters contains infinitely many clustering words if and only if, up to a permutation of letters, its directive word is $D^{(1)}D^{(2)}...D^{(r-1)}$, where $D^{(1)}$ is a finite word on the alphabet $\{a_1,a_2\}$, $D^{(2)}$ is a finite word on the alphabet $\{a_3,x_3\}$ with $x_3=a_1$ or $x_3=a_2$,  $D^{(3)}$ is a finite word on the alphabet $\{a_4,x_4\}$ with $x_4=a_3$ or $x_4=x_3$, ...,  $D^{(r-2)}$ is a finite word on the alphabet $\{a_{r-1},x_{r-1}\}$ with $x_{r-1}=a_{r-2}$ or $x_{r-1}=x_{r-2}$, $D^{(r-1)}$ is a one-sided infinite word on the alphabet $\{a_r,x_r\}$ with $x_r=a_{r-1}$ or $x_r=x_{r-1}$. This can be proved in the same way as Theorem \ref{thepi}.


\begin{thebibliography}{40}



\bibitem{ar} P. ARNOUX, G. RAUZY: Repr\'esentation g\'eom\'etrique de
suites
de complexit\'e $2n+1$ (French), {\it Bull. Soc. Math. France}  119 (1991),
p. 199--215.

\bibitem{be} J. BERSTEL: Sturmian and episturmian words (a survey of some recent results), in: Algebraic informatics, p. 23--47,  {\it Lecture Notes in Comput. Sci.} 4728, Springer, Berlin, 2007. 

\bibitem{bal} S. BRLEK, A. FROSINI, I. MANCINI, E. PERGOLA, S. RINALDI:  Burrows-Wheeler transform of words defined by morphisms, in: Combinatorial Algorithms, IWOCA 2019, p. 393--404, {\it Lecture Notes in Computer Science} 11638, Springer, Cham, 2019. 

\bibitem{brl} S. BRLEK, S. LI: On the number of squares in a finite word,  arXiv:2204.10204.

\bibitem{bw} M. BURROWS, D.J. WHEELER: A block-sorting lossless data compression algorithm, {\it Technical Report 124} (1994), Digital Equipment Corporation.

\bibitem{cc} J. CASSAIGNE, N. CHEKHOVA: Fonctions de récurrence des suites d'Arnoux-Rauzy et réponse à une question de Morse et Hedlund, (French) [Recurrence functions of Arnoux-Rauzy sequences and answer to a question posed by Morse and Hedlund], in: Numération, pavages, substitutions, {\it Ann. Inst. Fourier (Grenoble)} 56 (2006), p. 2249--2270. 

\bibitem{cfme} J. CASSAIGNE, S. FERENCZI, A. MESSAOUDI: Weak mixing and eigenvalues for Arnoux-Rauzy sequences, {\it Ann. Inst. Fourier (Grenoble)} 58 (2008),  p. 1983--2005.

\bibitem{zdl} A. DE LUCA, M. EDSON, L.Q.ZAMBONI: Extremal values of semi-regular continuants and codings of interval exchange transformations, {\it Mathematika} 69 (2023), p. 432-457.

\bibitem{djp} X. DROUBAY, J. JUSTIN, G. PIRILLO,  Episturmian words and some constructions of de Luca and Rauzy, {\it Theoret. Comput. Sci.} 255 (2001), p.  539--553. 


  \bibitem{fie} S. FERENCZI: A generalization of the self-dual induction to every interval exchange transformation, {\it Ann. Inst. Fourier (Grenoble)} 64 (2014), p. 1947--2002.
  
 \bibitem{fhuz} S. FERENCZI, P. HUBERT, L.Q. ZAMBONI: Languages of general interval exchange transformations, arXiv: 2212.01024.

\bibitem{fz1} S. FERENCZI, L.Q. ZAMBONI: Structure of K-interval exchange transformations: induction, trajectories, and distance theorems, {\it J. Anal. Math.} 112 (2010), p.  289--328.

\bibitem{fz3} S. FERENCZI, L.Q. ZAMBONI: Languages of k-interval exchange transformations, {\it Bull. Lond. Math. Soc.} 40 (2008), p. 705--714.

\bibitem{fz4} S. FERENCZI, L.Q. ZAMBONI: Clustering words and interval exchanges, {\it J. Integer Seq.} 16 (2013), Article 13.2.1, 9 pp.

\bibitem{gjp}  A. GLEN, J. JUSTIN: Episturmian words: a survey, {\it Theor. Inform. Appl.}  43 (2009), p. 403--442.

\bibitem{jz} O. JENKINSON, L. Q. ZAMBONI: Characterisations of balanced words via orderings, {\it Theoret. Comput. Sci.} 310 (2004), p. 247--271.

\bibitem{mrs} S. MANTACI, A. RESTIVO, M. SCIORTINO:  Burrows-Wheeler transform and Sturmian words, {\it Inform. Process. Lett.} 86 (2003), p.241--246. 

\bibitem{hm} M. MORSE, G.A. HEDLUND: Symbolic dynamics II. Sturmian trajectories, {\it  Amer. J. Math.} 62 (1940), p. 1--42.

\bibitem{sp} J. SIMPSON, S. J. PUGLISI: Words with simple Burrows-Wheeler Transforms, {\it The Electronic Journal of Combinatorics} 15 (2008), Research Paper 83, 17pp.





\end{thebibliography}
\end{document}